\documentclass{article}

\usepackage[preprint]{neurips_2021}

\usepackage[utf8]{inputenc} 
\usepackage[T1]{fontenc}    
\usepackage{hyperref}       
\usepackage{url}            
\usepackage{booktabs}       
\usepackage{amsfonts}       
\usepackage{nicefrac}       
\usepackage{microtype}      
\usepackage{xcolor}         

\usepackage{algorithm}
\usepackage{algorithmic}
\usepackage{amsmath}
\usepackage{mathrsfs}
\usepackage{amssymb}
\usepackage{caption}
\usepackage{subfig}
\usepackage{graphicx}

\title{Iteration Complexity of an Infeasible Interior Point Methods for Seconder-order Cone Programming and its Warmstarting}

\author{%
  Yushu Chen \thanks{Also at National Supercomputing Center in Wuxi, Jiangsu, China. \texttt{chenyushu@mail.tsinghua.edu.cn}} \space \space and Guangwen Yang \thanks{Corresponding author. Also at National Supercomputing Center in Wuxi, Jiangsu, China. \texttt{ygw@mail.tsinghua.edu.cn}}\\
  Department of Computer Science and Technology\\
  Tsinghua University\\
  Beijing, China, 100084\\
  \AND
  Lu Wang, Qingzhong Gan, and Haipeng Chen \\
  Shanghai Aerospace Control Technology Institute \\
  Shanghai, China, 201108\\
}

\begin{document}
\bibliographystyle{plainnat}

\maketitle
\begin{abstract}
  This paper studies the worst case iteration complexity of an infeasible interior point method (IPM) for seconder order cone programming (SOCP), which is more convenient for warmstarting compared with feasible IPMs. The method  studied bases on the homogeneous and self-dual model and the Monteiro-Zhang family of searching directions. Its worst case iteration complexity is $O\left(k^{1/2}\log\left(\epsilon^{-1}\right)\right)$, to reduce the primal residual, dual residual, and complementarity gap by a factor of $\epsilon$, where $k$ is the number of cone constraints. The result is the same as the best known result for feasible IPMs. The condition under which warmstarting improves the complexity bound is also studied.
\end{abstract}

\section{Introduction}

Warm starting of IPMs is widely perceived to be difficult \citep{Potra2000}. When the initial value is close to the boundary and not well-centered, the IPMs usually convergences very slow. Compared with feasible IPMs, since infeasible IPMs can be initialized trivially under the cone constraints without satisfying the equity constrains, they are more convenient for warm starting \citep{Skajaa2013}. In order to study the convergence of warm starting schemes for SOCP, a better understanding on the iteration complexity of infeasible IPMs is required. This paper studies the iteration-complexity of a infeasible interior point method for SOCP, and present the condition under which warmstarting improves the complexity.

Most of the existing works on the worst-case iteration complexity of IPMs for SOCP focus on feasible IPMs, requiring that the iterations satisfy the equity constraints. Among these works, \citet{Nesterov1997, Nesterov1998} present the Nesterov-Todd (NT) searching direction and show that the short-step path-following algorithm \citep{Kojima1989} using the direction has $O\left(k^{1/2}\log\left(\epsilon^{-1}\right)\right)$  iteration complexity, where $k$ is the number of cones. \citet{Tsuchiya1999} further studies the convergence of path-following algorithms for SOCP, and shows that: (1) using the HRVW/KSH/M directions \citep{Helmberg1996, Kojima1997, Monteiro1997}, the iteration complexity of the shortstep, semi-long step and long step methods are  $O\left(k^{1/2}\log\left(\epsilon^{-1}\right)\right)$, $O\left(k\log\left(\epsilon^{-1}\right)\right)$, and $O\left(k^{3/2}\log\left(\epsilon^{-1}\right)\right)$, respectively; (2) using the NT directions, the iteration complexity of the methods are  $O\left(k^{1/2}\log\left(\epsilon^{-1}\right)\right)$, $O\left(k\log\left(\epsilon^{-1}\right)\right)$, and $O\left(k\log\left(\epsilon^{-1}\right)\right)$ iterations, respectively. \citet{Monteiro2000} extend the $O\left(k^{1/2}\log\left(\epsilon^{-1}\right)\right)$ complexity bound to the short-step method and the Mizuno-Todd-Ye predictor-corrector method using the Monteiro-Zhang (MZ) family \citep{Monteiro1997, Monteiro1998, Zhang1998} family of directions, that includes the HRVW/KSH/M directions and NT directions as special cases. \citet{Mohammad2018} show a Newton method converges quadraticly, when initialized using a IPM solution with a sufficiently small complementary gap. 

Compared with feasible IPMs, iteration complexity of infeasible IPMs for SOCP are not thoroughly investigated. \citet{Bharath2006} study the convergence of infeasible IPMs for conic programing that includes SOCP as a special case. They show using the $\mathcal{N}_{-\infty}$ central path neighborhood, the method using the NT directions takes $O\left(k^2\log\left(\epsilon^{-1}\right)\right)$ iterations, and the method using th XS or SX directions \citep{Zhang1998} takes $O\left(k^{5/2}\log\left(\epsilon^{-1}\right)\right)$ iterations.

In this paper, we show an infeasible IPM for SOCP has $O\left(k^{1/2}\log\left(\epsilon^{-1}\right)\right)$ iteration complexity, that is the same as the best known result of feasible IPMs. We also present the condition under which warmstarting improves the bound. The method studied is an extension of the well know short-step path-following algorithm for linear programming, carried over to the context of SOCP. The modifications are as follows: firstly, the iterates are not required to be feasible to the equity constraints, but only be in the interior of the cone constraints, so that it is convenient for warmstarting; secondly, a homogeneous and self-dual (HSD) model \citep{Nemirovskii1993, Peng2003} is applied, therefore an infeasibility certificate can be obtained when either the primal or dual problem is infeasible; thirdly, the MZ family of searching directions is employed for acceleration. The warmstarting scheme analyzed is presented in a companion paper \citep{Chen2022}, which uses initial values modified from inexact solutions of previous problems to save computation.

The paper is organized as follows. Section \uppercase\expandafter{\romannumeral2} introduces the SOCP problem and the IPM for solving it. Section \uppercase\expandafter{\romannumeral3} studies the iteration complexity for the IPM to reduce the primal residual, dual residual, and complementarity gap. Section \uppercase\expandafter{\romannumeral4} studies the conditions that warm starting can improve the worst case complexity compared with cold starting.

\subsection{Notation and terminology}
The following notations are used throughout the paper. $\mathcal{R}_+$ denotes the set of positive real numbers, $\mathcal{R}^n$ denotes the $n$-dimensional Euclidean space, and $\mathcal{R}^{m\times n}$ denotes the set of all $m\times n$ matrices with real entries. $\mathbf{I}_n$ denotes the identity matrix of order $n$. $\mathbf{G}\succ 0$ means $\mathbf{G}$ is symmetric and positive definite, and $\mathbf{G}\succeq 0$ means $\mathbf{G}$ is symmetric and positive semi-definite. For a square matrix $\mathbf{M}$ with all real eigenvalues, its largest and smallest eigenvalues are denoted by $\lambda_{\max}(\mathbf{M})$ and $\lambda_{\min}(\mathbf{M})$. For a set $\mathcal{K}$, $\text{int}(\mathcal{K})$ denotes its interior. $\text{diag}(\cdots)$ denotes a diagonal matrix using a given vector as diagonal, and $\text{blkdiag}(\cdots)$ means a block diagonal matrix with given matrix blocks. $\left\|\cdots \right\|$ without specification is used instead of $\left\|\cdots \right\|_2$ for simplicity. With a little abuse of notations, when concatenating column vectors $\mathbf{x}$ and $\mathbf{y}$, we use $(\mathbf{x},\mathbf{y})$ instead of $(\mathbf{x}^T,\mathbf{y}^T)^T$ for simplicity.

\section{The SOCP problem and the infeasible IPM}

In this section, we introduce the seconder-order cone programming (SOCP) problem, some important concepts about the problem, and the infeasible IPM for SOCP studied in this paper.

\subsection{The SOCP problem}
SOCP minimize a linear function over the intersection of an affine set and the direct product of linear cones (LC) and second-order cones (SOC). 

The linear cones $\mathcal{K}_{L}^{n}$ and second-order cones $\mathcal{K}_{S}^{n}$ are defined as
\begin{equation}
\begin{aligned}
&\mathcal{K}_{L}^{n}=\{\mathbf{v }\in {{\mathbf{R}}^{n}}:\mathbf{v }\ge 0\}, n\ge 1, \\
&\mathcal{K}_{S}^{n}=\{\mathbf{v }\in {{\mathbf{R}}^{n}}:{{\mathbf{v }}_{1}}\ge \left\| {{\mathbf{v }}_{2:n}} \right\|\}, n\ge2,
\label{LC_SOC}
\end{aligned}
\end{equation}
The 1-dimensional SOC is defined as $\mathcal{K}_{S}^{1}=\mathcal{K}_{L}^{1}$. 

Then, the standard form of SOCP is defined as
\begin{equation}
\begin{aligned}
&\min_{\mathbf{x}}\mathbf{c}^T\mathbf{x},\\
&\text{s.t.}\\
&\mathbf{A}\mathbf{x}=\mathbf{b},\\
&\mathbf{x}\in\mathcal{K},\\
&\mathcal{K} =\mathcal{K}_{L}^{l}\times \mathcal{K}_{S}^{{{n}_{l+1}}}\times \mathcal{K}_{S}^{{{n}_{l+2}}}\cdots \times \mathcal{K}_{S}^{{{n}_{l+m}}},
\label{SOCP}
\end{aligned}
\end{equation}
where $\mathbf{A} \in \mathcal{R}^{p\times n}$, $\text{rank}(\mathbf{A})=p\le n$. Without loss of generality, we assume that in the solution variable the linear cone is arranged at first, followed with the SOCs. The number of cones is
\begin{equation}
k=l+m,
\label{NUM_CON}
\end{equation}
where the $l$-dimensional LC $\mathcal{K}_{L}^{n}$ is viewed as $l$ 1-dimensional LCs, so that $n_i=1$ for $i\in \left\{1,\cdots,l\right\}$. The $i$-th cone is also mentioned as $\mathcal{K}^i$ for simplicity, and correspondingly  $\mathcal{K}=\mathcal{K}^1\times\cdots\times\mathcal{K}^k$. For a vector $\mathbf{v}\in \mathcal{K}$, we use $\mathbf{v}^{(i)}$ to denote the subvector in the cone $\mathcal{K}^i$. For $\mathbf{v}\in \mathcal{R}^n-\mathcal{K}$, $\mathbf{v}^{(i)}$ also denotes the subvector corresponding to the position of $\mathcal{K}^i$.

The dual problem of SOCP Eq. (\ref{SOCP}) is
\begin{equation}
\begin{aligned}
&\max_{\mathbf{y}}\mathbf{b}^T\mathbf{y}\\
&\text{s.t.}\\
&\mathbf{A}^T\mathbf{y}+\mathbf{s}=\mathbf{c},\\
&\mathbf{s}\in\mathcal{K}.\\
\label{SOCP_DUAL}
\end{aligned}
\end{equation}

\subsection{Important concepts about SOCP}

This subsection introduce some important concepts about SOCP for later use.

Firstly, we introduce the Euclidean Jordan algebra associated with the SOC (see e.g. \citep{Faraut1994,Faybusovich1997a,Faybusovich1997b}), which is crucial in the IPMs for SOCP. The algebra for the SOC $\mathcal{K}^i$ is defined as
\begin{equation}
	\mathbf{x}^{(i)} \circ \mathbf{s}^{(i)} =\left({\left(\mathbf{x}^{(i)} \right)}^T \mathbf{s}^{(i)}, \mathbf{x}_1^{(i)}\mathbf{s}_{2:n_i}^{(i)}+\mathbf{s}_1^{(i)}\mathbf{x}_{2:n_i}^{(i)} \right), 
	\forall \mathbf{x}^{(i)},\mathbf{s}^{(i)}\in \mathcal{R}^{n_i}
	\label{DEF_CIRC_SOC}
\end{equation}
The Euclidean Jordan algebra for $\mathcal{K}$ is defined as
\begin{equation}
	\mathbf{x} \circ \mathbf{s} =\left(\mathbf{x}^{(1)}\circ\mathbf{s}^{(1)}, \cdots, \mathbf{x}^{(k)}\circ\mathbf{s}^{(k)}\right),
	\forall \mathbf{x},\mathbf{s}\in \mathcal{R}^{n}
	\label{DEF_CIRC_K}
\end{equation}

From now on, we assume that the space $\mathcal{R}^{n}$ is endowed with the Euclidean Jordan algebra. 

The unit element of the algebra is $\mathbf{e}=\left(\mathbf{e}^{(1)},\cdots,\mathbf{e}^{(k)}\right)$, where $\mathbf{e}^{(i)}=(1,0,\cdots,0)^T, i\in{1,\cdots,k}$, because $\mathbf{v}\circ\mathbf{e}=\mathbf{v}$. Then, the algebra can be written as multiplication of arrow head matrices and the unit element, as
\begin{equation}
	\mathbf{x} \circ \mathbf{s} = \text{mat}(\mathbf{x}) \text{mat}(\mathbf{s})\mathbf{e},
	\label{REL_CIRC_MAT}
\end{equation}
where the arrow head matrix for $\mathbf{v}\in\mathcal{R}^n$ is defined as
\begin{equation}
\begin{aligned}
&\text{mat}(\mathbf{v})=\text{blkdiag}\left(\text{mat}\left(\mathbf{v}^{(1)}\right)\cdots,\text{mat}\left(\mathbf{v}^{(k)}\right)\right),\\
&\text{mat}\left(\mathbf{v}^{(i)}\right)\triangleq \begin{pmatrix}
   {\mathbf{v}^{(i)}_{1}} & \left(\mathbf{v}^{(i)}_{2:n}\right)^{T}  \\
   \mathbf{v}^{(i)}_{2:n} & {\mathbf{v}^{(i)}_1}{{\mathbf{I}}_{n_i-1}} 
\end{pmatrix}, i\in{1,\cdots,k}.
\label{ARROW_MAT}
\end{aligned}
\end{equation}

It is easy to verify that the eigenvalues of the arrow head matrix satisfies
\begin{equation}
\lambda_{\min}(\text{mat}(\mathbf{v}))=\mathbf{v}_1-\|\mathbf{v}_{2:n}\|, \quad
\lambda_{\max}(\text{mat}(\mathbf{v}))=\mathbf{v}_1+\|\mathbf{v}_{2:n}\|.
\end{equation}
Consequently, we define the notation
\begin{equation}
\lambda_{\min}(\mathbf{v}) \triangleq \mathbf{v}_1-\|\mathbf{v}_{2:n}\|, \quad
\lambda_{\max}(\mathbf{v}) \triangleq \mathbf{v}_1+\|\mathbf{v}_{2:n}\|.
\end{equation}
For a vector $\mathbf{v}\in \mathcal{K}$, $\lambda_i^{1}(\mathbf{v})$ and $\lambda_i^{2}(\mathbf{v})$ are used as abbreviations of $\lambda_{\min}(\mathbf{v}^{(i)})$ and $\lambda_{\max}(\mathbf{v}^{(i)})$, where $\mathbf{v}^{(i)}$ denotes the part of $\mathbf{v}$ in the $i$-th cone. 

Secondly, we introduce the homogeneous and self-dual (HSD) model and the central path. The HSD model combines the  primal and dual problem, and enables to solved them simultaneously. The model is given by
\begin{equation}
\begin{aligned}
&\min \tilde{\beta} \nu \\						
&\text{s.t.}\\
&\mathbf{Ax}-\mathbf{b}\tau -\tilde{\mathbf{r}}_p\nu =0,\\				
&-{{\mathbf{A}}^{T}}\mathbf{y}+\mathbf{c}\tau -\mathbf{s}-\tilde{\mathbf{r}}_d\nu =0,\\			
&{{\mathbf{b}}^{T}}\mathbf{y}-{{\mathbf{c}}^{T}}\mathbf{x}-\kappa -\tilde{r}_g\nu =0,\\			
&\tilde{\mathbf{r}}_p^{T}\mathbf{y}+\mathbf{r}_{d,0}^{T}\mathbf{x}+\tilde{r}_g\tau =-\tilde{\beta},\\			
&(\mathbf{x},\mathbf{y},\mathbf{s},\kappa,\tau) \in \mathcal{C}, \mathcal{C}=\mathcal{K} \times \mathcal{R}^p \times \mathcal{K} \times \mathcal{K}_L^1 \times \mathcal{K}_L^1,
\label{HSD}
\end{aligned}
\end{equation}	
where $\nu \in {{R}_{+}}$ is free, and  ${{\mathbf{r}}_{p}}\in {{R}^{p}}$,  ${{\mathbf{r}}_{d}}\in {{R}^{n}}$, ${{r}_{g}}\in R$, $\beta \in R$ are residuals defined by
\begin{equation}
\begin{aligned}
&\tilde{\mathbf{r}}_p\triangleq \left( \mathbf{A}{{\mathbf{x}}_{0}}-\mathbf{b}{{\tau }_{0}} \right)/{{\nu }_{0}},\\
&\tilde{\mathbf{r}}_d\triangleq \left( -{{\mathbf{A}}^{T}}{{\mathbf{y}}_{0}}+\mathbf{c}{{\tau }_{0}}-{{\mathbf{s}}_{0}} \right)/{{\nu }_{0}},\\    			
&\tilde{r}_g\triangleq \left( {{\mathbf{b}}^{T}}{{\mathbf{y}}_{\mathbf{0}}}-{{\mathbf{c}}^{T}}{{\mathbf{x}}_{0}}-{{\kappa }_{0}} \right)/{{\nu }_{0}},\\ 				
&\tilde{\beta} \triangleq -\left( \mathbf{r}_{p}^{T}{{\mathbf{y}}_{0}}+\mathbf{r}_{d}^{T}{{\mathbf{x}}_{0}}+{{r}_{g}}{{\tau }_{0}} \right), 
\label{RESID}
\end{aligned}
\end{equation}	           
${{\nu }_{0}}$ is typically set as 1.

The solution of the HSD model fully describes the property of the primal problem and dual problem, that is summarized by the following theorem.

\textbf{Theorem 1}. Let $\left( \mathbf{x},\mathbf{y},\mathbf{s},\kappa ,\tau  \right)$ be an optimal solution of the HSD model Eq. (\ref{HSD}). Then, the following relations hold \citep{Terlaky2002,Wang2003}:

a) $\mathbf{x}^T\mathbf{s}=0$ and $\kappa\tau=0$;

b) If $\tau>0$, then $\left( \mathbf{x},\mathbf{y},\mathbf{s}\right)/\tau$ is the optimal solution of the primal problem Eq. (\ref{SOCP}) and the dual problem Eq. (\ref{SOCP_DUAL});

c) If $\kappa>0$, then at least one of the inequities $\mathbf{b}^T \mathbf{y} >0$ and $\mathbf{c}^T \mathbf{x} <0$ hold. If the first inequality holds then the primal problem Eq. (\ref{SOCP}) is infeasible. If the second inequality holds then the dual problem Eq. (\ref{SOCP_DUAL}) is infeasible.

\textit{Proof}. see corollary 1.3.1 of \citep{Wang2003}.

The central path is defined by 
\begin{equation}
\begin{aligned}
 & \mathbf{Ax}-\mathbf{b}\tau -{{\mathbf{r}}_{p}}\nu =0, \\ 
 & -{{\mathbf{A}}^{T}}\mathbf{y}+\mathbf{c}\tau -\mathbf{s}-{{\mathbf{r}}_{d}}\nu =0, \\ 
 & {{\mathbf{b}}^{T}}\mathbf{y}-{{\mathbf{c}}^{T}}\mathbf{x}-\kappa -{{r}_{g}}\nu =0, \\ 
 & \mathbf{x}\circ\mathbf{s}=\mu\mathbf{e},  \\ 
 & \kappa \tau =\mu
\label{CENTRAL_PATH}
\end{aligned}
\end{equation}
The central path approaches the optimal solution of the HSD model as $\mu \to 0$, that corresponds to the solution or an infeasibility certificate of the SOCP problem.

Thirdly, we introduce the scaling matrix, that can be applied to scale the searching directions. \citet{Monteiro2000} introduces a group of scaling auto-morphism that maps the cone $\mathcal{K}^i$ onto itself, as the following group of matrices
\begin{equation}
\mathcal{G}^i \triangleq 
\left\{\theta^i \mathbf{G}^i: \theta^i>0, \mathbf{G}^i \in \mathbf{R}^{n_i\times n_i}, 
(\mathbf{G}^i)^T \mathbf{Q}^i \mathbf{G}^i = \mathbf{Q}^i \right\},
\label{GROUP_AUTO_MORPHISM}
\end{equation}
where $\theta^i$ is a positive scaling value, and 
\begin{equation}
\mathbf{Q}^i = 
\begin{pmatrix}
1 & 0 \\
0 & -\mathbf{I}_{n_i-1}\\
\end{pmatrix}.
\label{DEF_Q}
\end{equation}

The set of scaling matrice for a vector in $\mathcal{K}$ is defined as
\begin{equation}
\mathcal{G} \triangleq 
\left\{
\Theta \mathbf{G}: \mathbf{G}=\text{blkdiag}\left( {{\mathbf{G}}^{1}},\cdots ,{{\mathbf{G}}^{k}}\right), 
\mathbf{G}^i \in \mathbf{R}^{n_i\times n_i}, (\mathbf{G}^i)^T \mathbf{Q}^i \mathbf{G}^i = \mathbf{Q}^i
\right\}
\end{equation}
where $i\in \left\{1,\cdots,k\right\}$ and
\begin{equation}
    \Theta =\text{blkdiag}\left({\theta }^{1}\mathbf{I}_{n_1},\cdots ,{\theta }^{k}\mathbf{I}_{n_k} \right).
\end{equation}

The solution variables $\mathbf{x}$ and $\mathbf{s}$ are scaled by
\begin{equation}
    \mathbf{\bar{x}}=(\Theta\mathbf{G})^T \mathbf{x}, \mathbf{\bar{s}}={{(\Theta \mathbf{G})}^{-1}}\mathbf{s},
    \label{DEF_SCAL_VAR}
\end{equation}
 For convenience, we denote a scaling matrix by
\begin{equation}
\mathbf{D}=(\mathbf{\Theta}\mathbf{G})^{-1}.
\label{DEF_MAT_D}
\end{equation}

The following proposition gives several important properties of the scaling matrix.

\textbf{Proposition 1} \citep{Andersen2003}. Let $\mathbf{D}\in \mathcal{G}$, then
\begin{equation}
\begin{aligned}
&(\mathbf{x}^{(i)})^T \mathbf{s}^{(i)} = (\mathbf{\bar{x}}^{(i)})^T \mathbf{\bar{s}}^{(i)},\\
&\theta_i^2 (\mathbf{x}^{(i)})^T \mathbf{Q}^i \mathbf{x}^{(i)} = (\mathbf{\bar{x}}^{(i)})^T \mathbf{Q}^i \mathbf{\bar{x}}^{(i)}\\
&\theta_i^{-2} (\mathbf{s}^{(i)})^T \mathbf{Q}^i \mathbf{s}^{(i)} = (\mathbf{\bar{s}}^{(i)})^T \mathbf{Q}^i \mathbf{\bar{s}}^{(i)},\\
&\mathbf{x} \in \mathcal{K} \Leftrightarrow \mathbf{\bar{x}} \in \mathcal{K}, \text{ and } 
\mathbf{x} \in \text{int} (\mathcal{K}) \Leftrightarrow \mathbf{\bar{x}} \in \text{int} (\mathcal{K}).
\label{PROP_SCAL_MAT}
\end{aligned}
\end{equation}
\textit{Proof}. See Lemma 3 of \citep{Andersen2003}.

For a vector $\mathbf{v}\in \text{int}\mathbf(K)$, \citet{Tsuchiya1999} introduces a special scaling matrix $\mathbf{T}_v$ defined as 
\begin{equation}
\mathbf{T}_v \triangleq \text{blkdiag}\left(\mathbf{T}_{v^{(1)}},\cdots,\mathbf{T}_{v^{(k)}}\right),
\label{DEF_T}
\end{equation}
where
\begin{equation}
\begin{aligned}
&\mathbf{T}_{u} \triangleq \begin{pmatrix}
&\mathbf{u}_{1}    &\mathbf{u}_{2:n}^T\\
&\mathbf{u}_{2:n}    &\beta_{u}\mathbf{I} + \frac{\mathbf{u}_{2:n}\mathbf{u}_{2:n}^T}{\beta_{u}+\mathbf{u}_{1} }
\end{pmatrix}, 
\beta_{u}=\sqrt{\mathbf{u}_{1}^2-\|\mathbf{u}_{2:n}\|^2}
\label{DEF_T_CON}
\end{aligned}
\end{equation}
for $\mathbf{u} \in \mathcal{R}^n$.

\textbf{Proposition 2} \citep{Tsuchiya1999}. For a vector $\mathbf{v}\in \text{int}\mathbf(K)$, $\mathbf{T}_v$ is the unique scaling matrix which maps $\mathbf{e}$ to $\mathbf{T}_v$. Moreover, $\mathbf{T}_v\succ 0$.

\textit{Proof}. See Proposition 2.1 of \citep{Tsuchiya1999}.

Finally, it is ready to define the scaling-invariant neighborhoods of the central path. For $(\mathbf{x}, \mathbf{s}, \kappa, \tau) \in \text{int}(\mathcal{K}) \times \text{int}(\mathcal{K}) \times R_+ \times R_+$, define the central path coefficient
\begin{equation}
\mu(\mathbf{x},\mathbf{s},\kappa,\tau)=(\mathbf{x}^T\mathbf{s}+\kappa\tau)/(k+1).
\label{MU}
\end{equation}
With a little abuse of notation, we also use $\mu$ instead of $\mu(\mathbf{x},\mathbf{s},\kappa,\tau)$ for simplicity.
The central path distances are defined as
\begin{equation}
\begin{aligned}
&d_2(\mathbf{x},\mathbf{s},\kappa,\tau)=\sqrt{2}\left\|\mathbf{w}_{xs\kappa\tau}-\mu\mathbf{e}\right\|,\\
&d_\infty(\mathbf{x},\mathbf{s},\kappa,\tau)=\max\left(d_\infty(\mathbf{x},\mathbf{s}), \left|\kappa\tau-\mu\right|\right),
\label{DIST_CEN_PATH}
\end{aligned}
\end{equation}
where
\begin{equation}
\begin{aligned}
&\mathbf{w}_{xs}=\mathbf{T}_x\mathbf{s}, 
\mathbf{w}_{xs\kappa\tau}=\left(\mathbf{w}_{xs}, \kappa\tau\right),\\
&d_\infty(\mathbf{x},\mathbf{s})=\max_{i=1,\cdots,k, j=0,1}\left(\left|(\mathbf{w}_{xs})_1^{(i)}+\|(\mathbf{w}_{xs})_{2:n_i}^{(i)}\|-\mu\right|,\left|(\mathbf{w}_{xs})_1^{(i)}-\|(\mathbf{w}_{xs})_{2:n_i}^{(i)}\|-\mu\right|\right),
\label{DEF_W}
\end{aligned}
\end{equation}
for $\mathbf{v} \in \mathcal{K}$.

The neighborhoods of the central path determined by the above distance are
\begin{equation}
\begin{aligned}
&\mathcal{N}_2(\gamma)=\left\{(\mathbf{x},\mathbf{y},\mathbf{s},\kappa,\tau) \in \text{int}(\mathcal{C}): d_2(\mathbf{x},\mathbf{s},\kappa,\tau) \le \gamma \mu(\mathbf{x},\mathbf{s},\kappa,\tau) \right\},\\
&\mathcal{N}_\infty(\gamma)=\left\{(\mathbf{x},\mathbf{y},\mathbf{s},\kappa,\tau) \in \text{int}(\mathcal{C}): d_\infty(\mathbf{x},\mathbf{s},\kappa,\tau) \le \gamma \mu(\mathbf{x},\mathbf{s},\kappa,\tau) \right\},
\label{NEIGHBOR_CEN_PATH}
\end{aligned}
\end{equation}
where the parameter $\gamma \in (0,1)$. For $(\mathbf{x},\mathbf{y},\mathbf{s},\kappa,\tau) \in \text{int}(\mathcal{C})$, it is easy to show $d_\infty(\mathbf{x},\mathbf{s},\kappa,\tau)\le d_2(\mathbf{x},\mathbf{s},\kappa,\tau)$, so that $\mathcal{N}_2(\gamma) \subseteq \mathcal{N}_\infty(\gamma)$.

The following proposition shows the invariance property of $\mathbf{w}_{xs}$.

\textbf{Proposition 3} \citet{Tsuchiya1999}. Suppose that $(\mathbf{x},\mathbf{s}) \in \text{int}(\mathcal{K})\times \text{int}(\mathcal{K})$. Let $\mathbf{w}\triangleq \mathbf{w}_{xs}$ and $\mathbf{\bar{w}}\triangleq \mathbf{w}_{\bar{x}\bar{s}}$, $i\in {1,\cdots,k}$. Then:

a) $\mathbf{\bar{w}}_0^{(i)}=\mathbf{w}_0^{(i)}$ and $\|\mathbf{\bar{w}}_1^{(i)}\|=\|\mathbf{w}_1^{(i)}\|$; 

b) $\lambda_{\min}(\mathbf{\bar{w}}_0^{(i)})=\lambda_{\min}(\mathbf{w}_0^{(i)})$ and $\lambda_{\max}(\mathbf{\bar{w}}_0^{(i)})=\lambda_{\max}(\mathbf{w}_0^{(i)})$;

c) $d_2(\mathbf{\bar{x}},\mathbf{\bar{s}})=d_2(\mathbf{x},\mathbf{s})$ and $d_\infty(\mathbf{\bar{x}},\mathbf{\bar{s}})=d_\infty(\mathbf{x},\mathbf{s})$.

\textit{Proof}. See Proposition 2.4 of \citet{Tsuchiya1999}.

Using Proposition 3, it is easy to verify that for a point $(\mathbf{x}, \mathbf{y}, \mathbf{s}, \kappa, \tau)$ on the central path defined by Eq. (\ref{CENTRAL_PATH}), the distances $d_2(\mathbf{x},\mathbf{s},\kappa,\tau)=d_\infty(\mathbf{x},\mathbf{s},\kappa,\tau)=0$.

\subsection{The infeasible IPM for SOCP}
In this subsection, we briefly introduce the infeasible IPM for SOCP studied in the paper.

The SOCP problem can be solved by loosely track the central path towards the optimal solution of the HSD model. In such a way, the updating direction in each iteration is computed by solving
\begin{equation}
\begin{aligned}
&\mathbf{A}\Delta \mathbf{x}-\mathbf{b}\Delta \tau =-\left( 1-\nu  \right)\left( \mathbf{Ax}-\mathbf{b}\tau  \right),\\	
&-{{\mathbf{A}}^{T}}\Delta \mathbf{y}+\mathbf{c} \Delta\tau -\Delta \mathbf{s}=-\left( 1-\nu  \right)\left( -{{\mathbf{A}}^{T}}\mathbf{y}+\mathbf{c}\tau -\mathbf{s} \right),\\			 
&{{\mathbf{b}}^{T}}\Delta \mathbf{y}-{{\mathbf{c}}^{T}}\Delta \mathbf{x}-\Delta \kappa =-\left( 1-\nu  \right)\left( {{\mathbf{b}}^{T}}\mathbf{y}-{{\mathbf{c}}^{T}}\mathbf{x}-\kappa  \right),\\			 
&\mathbf{X}\Delta \mathbf{s}+\mathbf{S}\Delta \mathbf{x}=\nu \mu \mathbf{e}-\mathbf{XSe},\\			 
&\kappa \Delta \tau +\tau \Delta \kappa =\nu \mu -\kappa \tau,
\label{LINSYS_ORI}
\end{aligned}
\end{equation}
where $\mathbf{X}=\text{mat}(\mathbf{x})$, $\mathbf{S}=\text{mat}(\mathbf{s})$. $\mu$ is $\mu(\mathbf{x},\mathbf{s},\kappa,\tau)$ defined in Eq. (\ref{MU}), which is proportional to the complementarity gap $\mathbf{x}^T\mathbf{s}+\kappa\tau$.

Instead of using the searching direction obtained by solving Eq (\ref{LINSYS_ORI}) directly, the MZ family of drections are applied to speed up the convergence. The directions are calculated by solving a scaled linear system, as 
\begin{equation}
	\begin{aligned}
		&\mathbf{A}\Delta \mathbf{x}-\mathbf{b}\Delta \tau =-\left( 1-\nu  \right)\left( \mathbf{Ax}-\mathbf{b}\tau  \right),\\	
		&-{{\mathbf{A}}^{T}}\Delta \mathbf{y}+\mathbf{c}\Delta\tau -\Delta \mathbf{s}=-\left( 1-\nu  \right)\left( -{{\mathbf{A}}^{T}}\mathbf{y}+\mathbf{c}\tau -\mathbf{s} \right),\\			 
		&{{\mathbf{b}}^{T}}\Delta \mathbf{y}-{{\mathbf{c}}^{T}}\Delta \mathbf{x}-\Delta \kappa =-\left( 1-\nu  \right)\left( {{\mathbf{b}}^{T}}\mathbf{y}-{{\mathbf{c}}^{T}}\mathbf{x}-\kappa  \right),\\			 
		&\bar{\mathbf{X}} \mathbf{D} \Delta \mathbf{s} + \bar{\mathbf{S}}\mathbf{D}^{-T} \Delta \mathbf{x}
		=\nu \mu \mathbf{e} - \bar{\mathbf{X}}\bar{\mathbf{S}}\mathbf{e},\\			 
		&\kappa \Delta \tau +\tau \Delta \kappa =\nu \mu -\kappa \tau,
		\label{LINSYS_SCAL}
	\end{aligned}
\end{equation}
where 
\begin{equation}
	\bar{\mathbf{X}}=\text{mat}(\bar{\mathbf{x}}),\bar{\mathbf{S}}=\text{mat}(\bar{\mathbf{s}}), \bar{\mathbf{x}}=\mathbf{D}^{-T}\mathbf{x},\bar{\mathbf{s}}=\mathbf{D}\mathbf{s}, \mathbf{D} \in \mathcal{G}.
\end{equation}
The scaled linear system Eq (\ref{LINSYS_SCAL}) is obtained by using the scale variable $\bar{\mathbf{x}},\bar{\mathbf{s}}$ instead of $\mathbf{x}$ and $\mathbf{s}$ in the HSD model Eq. (\ref{HSD}).

The MZ family of directions include the HRVW/KSH/M directions and NT directions as special cases. For example, in the NT scaling, the scaling matrix $\mathbf{D}$ is symmetric and satisfies $\mathbf{x}=\mathbf{D}^2\mathbf{s}$, and the multiplication with $\mathbf{D}$ or $\mathbf{D}^2$ can be calculated in $O(n)$ complexity.

Then, the infeasible IPM for SOCP is presented as Algorithm 1.

\begin{algorithm}[tb]
	\caption{An Infeasible IPM for SOCP}
	\label{Algorithm1}
	\textbf{Input}: Constant coefficients $\gamma, \delta \in(0,1)$, a small positive threshold $\epsilon\in(0,1)$; the SOCP problem parameters $\mathbf{A},\mathbf{b},\mathbf{c}, \mathcal{K}$; initial value $(\mathbf{x}_0,\mathbf{y}_0,\mathbf{s}_0,\kappa_0,\tau_0) \in \mathcal{N}_2(\gamma)$. \\
	\textbf{Output}: the solution $(\mathbf{x},\mathbf{y},\mathbf{s},\kappa,\tau)$ and final status.
	\begin{algorithmic}[1] 
		\STATE Set $\nu=1-\delta/\sqrt{2(k+1)}$;
		\STATE Set $(\mathbf{x},\mathbf{y},\mathbf{s},\kappa,\tau)=(\mathbf{x}_0,\mathbf{y}_0,\mathbf{s}_0,\kappa_0,\tau_0)$;
		\WHILE {the stop criteria is not satisfied}
		\STATE Set $\mu=(\mathbf{x}^T\mathbf{s}+\kappa\tau)/(k+1)$;
		\STATE Compute $(\Delta\mathbf{x},\Delta\mathbf{y},\Delta\mathbf{s},\Delta\kappa,\Delta\tau)$ by solving the linear system Eq. (\ref{LINSYS_SCAL});
		\STATE Set $(\mathbf{x},\mathbf{y},\mathbf{s},\kappa,\tau)=(\mathbf{x},\mathbf{y},\mathbf{s},\kappa,\tau)+(\Delta\mathbf{x},\Delta\mathbf{y},\Delta\mathbf{s},\Delta\kappa,\Delta\tau)$;
		\ENDWHILE
	\STATE Select the final status according to the stop criteria.
	\end{algorithmic}
\end{algorithm}

In Algorithm 1, the stop criteria is as follows:
\begin{equation}
\begin{aligned}
& \left\| \mathbf{r}_p(\mathbf{x},\mathbf{y},\mathbf{s},\kappa,\tau) \right\|\le \epsilon \left\| \mathbf{r}_p(\mathbf{x}_0,\mathbf{y}_0,\mathbf{s}_0,\kappa_0,\tau_0) \right\|, \\ 
& \left\| \mathbf{r}_d(\mathbf{x},\mathbf{y},\mathbf{s},\kappa,\tau) \right\|\le \epsilon \left\| \mathbf{r}_d(\mathbf{x}_0,\mathbf{y}_0,\mathbf{s}_0,\kappa_0,\tau_0) \right\|, \\ 
& \mu \left( \mathbf{x},\mathbf{s},\kappa ,\tau  \right)\le \epsilon \mu \left( {{\mathbf{x}}_{0}},{{\mathbf{s}}_{0}},{{\kappa }_{0}},{{\tau }_{0}} \right).
\label{STOP_CRIT}
\end{aligned}
\end{equation}
where the primal residual $\mathbf{r}_p$ and the dual residual $\mathbf{r}_d$ are defined as
\begin{equation}
\mathbf{r}_p(\mathbf{x},\mathbf{y},\mathbf{s},\kappa,\tau)=\mathbf{Ax}-\tau \mathbf{b},
\mathbf{r}_d(\mathbf{x},\mathbf{y},\mathbf{s},\kappa,\tau)={{\mathbf{A}}^{T}}\mathbf{y}+\mathbf{s}-\tau \mathbf{c}.	
\end{equation}

For convience, we also define
\begin{equation} r_g(\mathbf{x},\mathbf{y},\mathbf{s},\kappa,\tau)={{\mathbf{b}}^{T}}\mathbf{y}-{{\mathbf{c}}^{T}}\mathbf{x}-\kappa.
\end{equation}

The criteria indicates that the primal residual, dual residual, and complementarity gap are reduce by a certain ratio $\epsilon$. When the criteria is satisfied, the final status is decided using an approximate version of Theorem 1: if $\tau \ge \epsilon\max(1,\kappa)$, then the problem is solved; if $\tau < \epsilon\max(1,\kappa)$, the primal problem is infeasible if $\mathbf{b}^T \mathbf{y} >0$, and the dual problem is infeasible if $\mathbf{c}^T \mathbf{x} <0$. Some important issues in implementation such as the maximum iteration number limit and rounding errors are not considered in the theoretical analysis of this paper.

Algorithm 1 obtains the solution for the HSD model. If the problem is solved, the solution for the SOCP problem is $(\mathbf{x}/\tau,\mathbf{y}/\tau,\mathbf{s}/\tau)$. In oder to guarantee the convergence of Algorithm 1, the parameters $\gamma, \delta$ should be chosen to satisfy certain conditions, which will be presented in Proposition 4 in the next section.

\section{Analysis of the iteration complexity}

In order to prove the polynomial convergence in the worst case, it is crucial to show two points: (1) the updating keeps the solution in a central path neighborhood, which is a subset of the interior of cone constraints; (2) when the solution is bounded in the central path neighborhood, the primal residual, dual residual, and complementarity gap can be reduced for a certain ratio in each iteration.

This section is organized as follows: firstly, an equivalent form of the HSD model is proposed to simplify the following analysis; secondly, the well-definedness of the searching directions is proved; thirdly, the improvement of the primal residual, dual residual, and complementarity gap is investigated; fourthly, several important results are presented to bound the updating in the central path neighborhood for the unscaled search directions; fifthly, the polynomial convergence for the unscaled search directions is proved; finally, the proof is extend to the MZ family of search directions.

\subsection{An equivalent form of the HSD model}

In the HSD model, $\kappa$ and $\tau$ are variables in $\mathcal{K}_L^1$, but they act differently from other variables in LCs. Although we find it is theoretically not difficult to consider them as additional terms, the analysis in such a way would be long and tedious. In order to simplify the analysis, we propose a equivalent form of the HSD model in this subsection.

Consider the following change of variables
\begin{equation}
\begin{aligned}
\hat{\mathbf{A}}=[\mathbf{A},-\mathbf{b}],
\hat{\mathbf{C}}=
\begin{pmatrix}\begin{aligned}
&0, &-\mathbf{c}\\
&\mathbf{c}^T, &0\\
\end{aligned}\end{pmatrix},
\hat{\mathbf{x}}=
\begin{pmatrix}
\mathbf{x}\\
\tau\\
\end{pmatrix},
\hat{\mathbf{s}}=
\begin{pmatrix}
\mathbf{s}\\
\kappa\\
\end{pmatrix},
\hat{\mathbf{e}}=
\begin{pmatrix}
\mathbf{e}\\
1\\
\end{pmatrix},
\label{CH_VAR_HSD}
\end{aligned}
\end{equation}
where $\hat{\mathbf{x}}, \hat{\mathbf{s}}\in \mathcal{K} \times \mathcal{K}_L^1$ and $\mathbf{y}\in\mathcal{R}^p$.

It is easy to verify that the linear system Eq. (\ref{LINSYS_ORI}) can be written as
\begin{equation}
\begin{aligned}
&\hat{\mathbf{A}}\Delta \hat{\mathbf{x}} =-\left( 1-\nu  \right)\hat{\mathbf{A}}\hat{\mathbf{x}},\\	
&\hat{\mathbf{A}}^{T}\Delta\mathbf{y}+\hat{\mathbf{C}} \Delta\hat{\mathbf{x}} +\Delta \hat{\mathbf{s}}=
-(1-\nu) \left(\hat{\mathbf{A}}^{T}\mathbf{y}+\hat{\mathbf{C}} \hat{\mathbf{x}} + \hat{\mathbf{s}} \right),\\
&\hat{\mathbf{X}}\Delta \hat{\mathbf{s}}+\hat{\mathbf{S}}\Delta \hat{\mathbf{x}}=\nu \mu \hat{\mathbf{e}}-\hat{\mathbf{X}}\hat{\mathbf{S}}\hat{\mathbf{e}},
\label{LINSYS_ORI_SIMP}
\end{aligned}
\end{equation}
where $\hat{\mathbf{X}}=\text{mat}(\hat{\mathbf{x}})$, $\hat{\mathbf{S}}=\text{mat}(\hat{\mathbf{s}})$. It is worth noting that the second and third equations in Eq. (\ref{LINSYS_ORI}) are combined into the second equation in Eq. (\ref{LINSYS_ORI_SIMP}).

The scaled linear system Eq. (\ref{LINSYS_SCAL}) is written as
\begin{equation}
\begin{aligned}
&\hat{\mathbf{A}}\Delta \hat{\mathbf{x}} =-\left( 1-\nu  \right)\hat{\mathbf{A}}\hat{\mathbf{x}},\\	
&\hat{\mathbf{A}}^{T}\Delta\mathbf{y}+\hat{\mathbf{C}} \Delta\hat{\mathbf{x}} +\Delta \hat{\mathbf{s}}=
-(1-\nu) \left(\hat{\mathbf{A}}^{T}\mathbf{y}+\hat{\mathbf{C}} \hat{\mathbf{x}} + \hat{\mathbf{s}} \right),\\			 
&\text{mat}(\hat{\mathbf{D}}^{-T}\hat{\mathbf{x}}) \hat{\mathbf{D}} \Delta \hat{\mathbf{s}} + \text{mat}(\hat{\mathbf{D}}\hat{\mathbf{s}}) \hat{\mathbf{D}}^{-T} \Delta \hat{\mathbf{x}}
=\nu \mu \hat{\mathbf{e}} - (\hat{\mathbf{D}}^{-T}\hat{\mathbf{x}}) \circ (\hat{\mathbf{D}}\hat{\mathbf{s}}),
\label{LINSYS_SCAL_SIMP}
\end{aligned}
\end{equation}
where the scaling matrix
\begin{equation}
\hat{\mathbf{D}}=
\begin{pmatrix}
\begin{aligned}
&\mathbf{D} &0\\
&0 &1
\end{aligned}
\end{pmatrix}.
\end{equation}

The central path coefficient is also simplified as
\begin{equation}
\mu=(\mathbf{x}^T\mathbf{s}+\kappa\tau)/(k+1)=\hat{\mathbf{x}}^T\hat{\mathbf{s}}/(k+1)
\label{MU_SIMP}
\end{equation}

Because the two linear systems are identical, the IPM solving the HSD model by Eq. (\ref{LINSYS_SCAL}) is identical to solving the model by Eq. (\ref{LINSYS_SCAL_SIMP}). Without loss of generality, the last LC $\mathcal{K}_L^1$ can be viewed as a SOC $\mathcal{K}_S^1$. Then, with a little abuse of notation, Eq. (\ref{LINSYS_SCAL_SIMP}) is a special case of the following generalized linear system
\begin{equation}
\begin{aligned}
&\mathbf{A}\Delta \mathbf{x} =-\left( 1-\nu  \right)\mathbf{A}\mathbf{x},\\	
&\mathbf{A}^{T}\Delta\mathbf{y}+\mathbf{C} \Delta\mathbf{x} +\Delta \mathbf{s}=
-(1-\nu) \left(\mathbf{A}^{T}\mathbf{y}+\mathbf{C} \mathbf{x} + \mathbf{s} \right),\\			 
&\bar{\mathbf{X}} \mathbf{D} \Delta \mathbf{s} + \bar{\mathbf{S}} \mathbf{D}^{-T} \Delta \mathbf{x}
=\nu \mu \mathbf{e} - \bar{\mathbf{x}} \circ \bar{\mathbf{s}},
\label{LINSYS_SCAL_GEN}
\end{aligned}
\end{equation}
where the variables are redefined for simplicity, as
\begin{equation}
\begin{aligned}
&\mathbf{A} \in \mathcal{R}^{p\times n}, \text{rank}(\mathbf{A})=p\le n, 
\mathbf{C} \in \mathcal{R}^{n\times n}, \mathbf{C}=-\mathbf{C}^T\\
&\mathcal{K} =\mathcal{K}_{L}^{l} \times \mathcal{K}_{S}^{{{n}_{l+1}}}\times \mathcal{K}_{S}^{{{n}_{l+2}}}\cdots \times \mathcal{K}_{S}^{{{n}_{l+m}}}, k=l+m, \\
&\mathbf{x} \in \mathcal{K}, \mathbf{y} \in \mathcal{R}^p, \mathbf{s} \in \mathcal{K},
\mathcal{C}=\mathcal{K} \times \mathcal{R}^p \times \mathcal{K},\\
&\mathbf{D}\in \mathcal{G}, \bar{\mathbf{x}}=\mathbf{D}^{-T}\mathbf{x},\bar{\mathbf{s}}=\mathbf{D}\mathbf{s},
\bar{\mathbf{X}}=\text{mat}(\bar{\mathbf{x}}),\bar{\mathbf{S}}=\text{mat}(\bar{\mathbf{s}}).
\label{DEF_VAR_GEN}
\end{aligned}
\end{equation}

The central path coefficient for the  generalized linear system Eq. (\ref{LINSYS_SCAL_GEN})  is defined as
\begin{equation}
\mu=\mu(\mathbf{x},\mathbf{s})=\mathbf{x}^T\mathbf{s}/k,
\label{MU_GEN}
\end{equation}
where $\mathbf{x}^T\mathbf{s}$ is the complementarity gap . The residuals are defined as
\begin{equation}
\hat{\mathbf{r}}_p=\mathbf{A}\mathbf{x},\hat{\mathbf{r}}_d=\mathbf{A}^{T}\mathbf{y}+\mathbf{C} \mathbf{x} + \mathbf{s}.
\label{RESIDUAL_GEN}
\end{equation}
For the special case Eq. (\ref{LINSYS_SCAL_SIMP}),  the generalized linear system Eq. (\ref{LINSYS_SCAL_GEN}) has 1 additional dimensional in both $\mathbf{x}$ and $\mathbf{s}$ compared with the original form Eq. (\ref{LINSYS_SCAL}), to substitute $\kappa$ and $\tau$. Meanwhile, its cone constraint $\mathcal{K}$ is also 1-dimensional higher. The scaling matrix $\mathbf{D}\in \mathcal{G}$ generates the MZ family of searching directions. The residual $\hat{\mathbf{r}}_p$ is identical to the primal residual $\mathbf{r}_p$ of the HSD model defined in Eq. (\ref{RESID}), $\hat{\mathbf{r}}_d$ is a combination of the dual residual and gap as $(\mathbf{r}_d,r_g)$, and the central path coefficient $\mu(\mathbf{x},\mathbf{s})$ is identical to $\mu(\mathbf{x},\mathbf{s},\kappa,\tau)$ defined in Eq. (\ref{MU}).

In the following, it sufficient to study the convergence properties to reduce the residuals $\hat{\mathbf{r}}_p$, $\hat{\mathbf{r}}_d$, and the complementarity gap which is proportional to $\mu$ by solving the linear system  Eq. (\ref{LINSYS_SCAL_GEN}), and the results apply to the process by solving Eq. (\ref{LINSYS_SCAL}).

For the generalized linear system Eq. (\ref{LINSYS_SCAL_GEN}), the central path distances are
\begin{equation}
\begin{aligned}
&d_2(\mathbf{x},\mathbf{s})=\sqrt{2}\left\|\mathbf{w}_{xs}-\mu\mathbf{e}\right\|,\\
&d_\infty(\mathbf{x},\mathbf{s})=\max_{i=1,\cdots,k, j=0,1}\left(\left|(\mathbf{w}_{xs})_0^{(i)}+\|(\mathbf{w}_{xs})_1^{(i)}\|-\mu\right|,\left|(\mathbf{w}_{xs})_0^{(i)}-\|(\mathbf{w}_{xs})_1^{(i)}\|-\mu\right|\right).
\label{DEF_DIST_CEN_PATH_GEN}
\end{aligned}
\end{equation}
The central path neighborhoods are
\begin{equation}
\begin{aligned}
&\mathcal{N}_2(\gamma)=\left\{(\mathbf{x},\mathbf{y},\mathbf{s}) \in \text{int}(\mathcal{C}): d_2(\mathbf{x},\mathbf{s}) \le \gamma \mu(\mathbf{x},\mathbf{s}) \right\},\\
&\mathcal{N}_\infty(\gamma)=\left\{(\mathbf{x},\mathbf{y},\mathbf{s}) \in \text{int}(\mathcal{C}): d_\infty(\mathbf{x},\mathbf{s}) \le \gamma \mu(\mathbf{x},\mathbf{s}) \right\},
\label{DEF_NEIGHBOR_CEN_PATH_GEN}
\end{aligned}
\end{equation}
where the parameter $\gamma \in (0,1)$. It is easy to verify the central path distances and neighborhoods are identical to the definition Eq. (\ref{DIST_CEN_PATH}) and Eq. (\ref{NEIGHBOR_CEN_PATH}) for the special case Eq. (\ref{LINSYS_SCAL_SIMP}). We also have  $d_\infty(\mathbf{x},\mathbf{s})\le d_2(\mathbf{x},\mathbf{s})$, and $\mathcal{N}_2(\gamma) \subseteq  \mathcal{N}_\infty(\gamma)$.

\subsection{The well-definedness of the searching directions}

\citet{Tsuchiya1999} and \citet{Monteiro2000} provide valuable tools for our study. For the convenience of analysis, we introduce Lemma 1-4 in \citep{Monteiro2000} as Lemma 1-4 (see the paper for the proofs).

Define the following notations as
\begin{equation}
	\mathbf{W}_{xs}  \triangleq   \text{mat} (\mathbf{w}_{xs}),
	\mathbf{R}(\mathbf{x},\mathbf{s})  \triangleq  \mathbf{T}_x \mathbf{X}^{-1} \mathbf{S} \mathbf{T}_x.
	\label{DEF_W_R}
\end{equation}

\textbf{Lemma 1} \citep{Monteiro2000}. For any $\mathbf{x}\in \text{int}(\mathcal{K})$, the matrice $\mathbf{X}$ and $\mathbf{T}_x$ satisfy:

a) 
\begin{equation}
\mathbf{X}-\mathbf{T}_x = \mathbf{U}_{x} =\text{blkdiag}\left(\mathbf{U}_{x^{(i)}}:i=1,\cdots,k \right),
\label{DEF_U}
\end{equation}
where 
\begin{equation}
	\mathbf{U}_{v}=
	\begin{pmatrix}
		&0	&0\\
		&0	&(\mathbf{v}_{1}-\beta_{v})\mathbf{P}_{v}
	\end{pmatrix}
\end{equation}
for $\mathbf{v}\in \mathcal{R}^n$, and $\mathbf{P}_{v}$ is the orthogonal projection matrix onto the subspace orthogonal to $\mathbf{v}_{2:n}$, as
\begin{equation}
	\mathbf{P}_{v}=\mathbf{I}_{n-1}-\mathbf{v}_{2:n}(\mathbf{v}_{2:n})^T / \left\|\mathbf{v}_{2:n}\right\|^2.
	\label{DEF_P}
\end{equation}
We also use the notation
\begin{equation}
	\mathbf{P}_{x}=\text{blkdiag}\left(
		\begin{pmatrix}
		&0	&0\\
		&0	&\mathbf{P}_{x^{(i)}}
		\end{pmatrix}
 	:i=1,\cdots,k\right).
	\label{DEF_PX}
\end{equation}

b) 
\begin{equation}
\begin{aligned}
&\mathbf{T}_x\mathbf{X}^{-1}=\mathbf{X}^{-1}\mathbf{T}_x
=\text{blkdiag}\left(\mathbf{I}-\mathbf{U}_{x^{(i)}}/\mathbf{x}^{(i)}_1: i=  1,\cdots,k\right),\\
&\mathbf{T}_x\mathbf{X}^{-1}\mathbf{e}=\mathbf{e};
\end{aligned}
\end{equation}

c) $\mathbf{X}$ and $\mathbf{T}_x$ commute and $\mathbf{X} \succeq \mathbf{T}_x$.

\textbf{Lemma 2} \citep{Monteiro2000}. $\mathbf{R}_{xs}$ satisfies
\begin{equation}
	\mathbf{R}_{xs} =\text{blkdiag}\left(\mathbf{R}_{xs}^i:i=1,\cdots,k \right),
\end{equation}
where
\begin{equation}
\mathbf{R}_{xs}^i \triangleq
\begin{pmatrix}
	&\dot{\mathbf{w}}_1			&\dot{\mathbf{w}}_{2:n_i}^T\\
	&\dot{\mathbf{w}}_{2:n_i}	&\tilde{\mathbf{R}}_{xs}^i
\end{pmatrix}
\end{equation}
with $\dot{\mathbf{w}}$ denotes $\mathbf{w}_{xs}^{(i)}$ and
\begin{equation}
\tilde{\mathbf{R}}_{xs}^i \triangleq \frac{1}{\mathbf{x}^{(i)}_1}
\left(\dot{\mathbf{w}}_{2:n_i}(\mathbf{x}^{(i)}_{2:n_i})^T+\beta_{x^{(i)}}^2\mathbf{s}^{(i)}_1 \mathbf{I}\right)
=\frac{\dot{\mathbf{w}}_{2:n_i}(\mathbf{x}^{(i)}_{2:n_i})^T}{\mathbf{x}^{(i)}_1}
+\left(\dot{\mathbf{w}}_1-\frac{\dot{\mathbf{w}}_{2:n_i}^T \mathbf{x}^{(i)}_{2:n_i}}{\mathbf{x}^{(i)}_1}\right)
\mathbf{I}.
\end{equation}
.

\textbf{Lemma 3} \citep{Monteiro2000}.  Let $(\mathbf{x},\mathbf{s}) \in \text{int}(\mathcal{K}) \times \text{int}(\mathcal{K})$ satisfy 
\begin{equation}
\max_{i,j} \left|\lambda_i^j(\mathbf{w}_{xs})-\dot{\nu}\right| \le \gamma \dot{\nu}
\label{LEMMA3_COND}
\end{equation}
 for some scalars $\gamma>0$ and $\dot{\nu}>0$. Then
\begin{equation}
\begin{aligned}
&\left\|\mathbf{R}_{xs}-\mathbf{W_{xs}}\right\| \le 2\gamma\dot{\nu},\\
&\left\|\mathbf{W}_{xs}-\dot{\nu} \mathbf{I}\right\| \le \gamma\dot{\nu},\\
&\left\|\mathbf{R}_{xs}-\dot{\nu} \mathbf{I}\right\| \le 3\gamma\dot{\nu},
\label{LEMMA3}
\end{aligned}
\end{equation}
\textbf{Lemma 4} \citep{Monteiro2000}. Let $(\mathbf{x},\mathbf{s}) \in \text{int}(\mathcal{K}) \times \text{int}(\mathcal{K})$ satisfy 
\begin{equation}
\left\|\mathbf{R}_{xs}-\dot{\nu} \mathbf{I}\right\| \le \dot{\tau} \dot{\nu},
\label{LEMMA4_COND1}
\end{equation}
for some scalars $\dot{\tau} \in (0,1)$ and $\dot{\nu}>0$. Assume that $(\mathbf{u},\mathbf{v})\in \mathcal{K} \times \mathcal{K}$ and $\mathbf{h} \in \mathcal{K}$ satisfy
\begin{equation}
\mathbf{S}\mathbf{u}+\mathbf{X}\mathbf{v}=\mathbf{h}, \mathbf{u}^T\mathbf{v} \ge 0,
\label{LEMMA4_COND2}
\end{equation}
and define $\delta_u \triangleq \left\|\mathbf{T}_x^{-1}\mathbf{u}\right\|$ and $\delta_v \triangleq \left\|\mathbf{T}_x \mathbf{v}\right\|$. Then,
\begin{equation}
\delta_u \le \frac{\left\|\mathbf{T}_x\mathbf{X}^{-1}\mathbf{h}\right\|}{(1-\dot{\tau})\nu},
\delta_v \le \frac{2\left\|\mathbf{T}_x\mathbf{X}^{-1}\mathbf{h}\right\|}{1-\dot{\tau}}
\label{LEMMA4}
\end{equation}

Using the above lemmas, we obtain the following result that shows the well-definedness of the search directions.

\textbf{Theorem 2}. Let $(\mathbf{x},\mathbf{y},\mathbf{s}) \in \text{int}(\mathcal{C})$ satisfy 
\begin{equation}
\max_{i,j}\left\|\lambda_i^j(\mathbf{w}_{xs})-\nu\right\| \le \gamma \nu
\label{THEOREM2_COND}
\end{equation}
 for $\gamma\in(0,1/3)$, and $
\mathbf{D} \in \mathcal{G}$. Then, the linear systems Eq. (\ref{LINSYS_SCAL_GEN}) has unique solution.

\textit{Proof}. In order to prove the linear system Eq. (\ref{LINSYS_SCAL_GEN}) has unique solution, it is equivalent to show the solution of the homogeneous system associated with it is 0. The homogeneous system is
\begin{equation}
\begin{aligned}
&\mathbf{A}\mathbf{d}_{hx} =0,\\		
&{{\mathbf{A}}^{T}} \mathbf{d}_{hy}+\mathbf{C} d_{hx} +\mathbf{d}_{hs}=0,\\
&\bar{\mathbf{X}}\mathbf{D}\mathbf{d}_{hs}+\bar{\mathbf{S}}\mathbf{D}^{-T}\mathbf{d}_{hx}=0,
\label{LINSYS_SCAL_GEN_H}
\end{aligned}
\end{equation}
where the solution is denoted by $(\mathbf{d}_{hx},\mathbf{d}_{hy},\mathbf{d}_{hs})$.

Because $\mathbf{D}\in\mathcal{G}$,  it is invertible. Since $\mathbf{C}=-\mathbf{C}^T$ by definition Eq. (\ref{DEF_VAR_GEN}), multiplying $\mathbf{d}_{hy}^T$ and $\mathbf{d}_{hx}^T$ to the first and second equations in Eq. (\ref{LINSYS_SCAL_GEN_H}) respectively, we can obtain
\begin{equation}
(\mathbf{D}^{-T}\mathbf{d}_{hx})^{T}(\mathbf{D}\mathbf{d}_{hs})=\mathbf{d}_{hx}^{T}\mathbf{d}_{hs}=0.
\end{equation}
By proposition 3 b), 
\begin{equation}
\max_{i,j}\left\|\lambda_i^j(\mathbf{w}_{\bar{x}\bar{s}})-\nu\right\| =\max_{i,j}\left\|\lambda_i^j(\mathbf{w}_{xs})-\nu\right\| \le \gamma. 
\end{equation}
By Lemma 3, Eq. (\ref{LEMMA4_COND1}) holds for $\dot{\tau}=3\gamma < 1$. By the last equation in Eq. (\ref{LINSYS_SCAL_GEN_H}), using $\mathbf{D}\mathbf{d}_{hs}$ and $\mathbf{D}^{-1}\mathbf{d}_{hx}$ as $\mathbf{u}$ and $\mathbf{v}$ in Lemma 4 with $\textbf{h}=0$, we obtain 
\begin{equation}
    \mathbf{d}_{hx}=0,\mathbf{d}_{hs}=0.
\end{equation}
Then, by the second equation in Eq. (\ref{LINSYS_SCAL_GEN_H}),
\begin{equation}
    {{\mathbf{A}}^{T}} \mathbf{d}_{hy}=0.
\end{equation}
Since $\textbf{A}$ has full row rank, we obtain
\begin{equation}
\mathbf{d}_{hy}=0.
\end{equation}
Consequently, the linear system Eq. (\ref{LINSYS_SCAL_GEN_H}) has unique solution. Then, the linear system Eq. (\ref{LINSYS_SCAL_GEN}) also has unique solution. \textbf{Q.E.D}.

\subsection{The improvements in each iteration}

Now we investigate the improvement of the primal residual, dual residual and complementarity gap in each iteration. For a step size $\alpha$, define
\begin{equation}
\begin{aligned}
&\mathbf{x}(\alpha) \triangleq \mathbf{x} + \alpha \Delta \mathbf{x},
\mathbf{y}(\alpha) \triangleq \mathbf{y} + \alpha \Delta \mathbf{y},
\mathbf{s}(\alpha) \triangleq \mathbf{s} + \alpha \Delta \mathbf{s}.
\label{INCRE_ALPHA}
\end{aligned}
\end{equation}

The following theorem states the improvement in each iteration to solve the linear system  Eq. (\ref{LINSYS_SCAL_GEN}). It also shows that the updating directions of $\mathbf{x}$ and $\mathbf{s}$ are orthogonal.

\textbf{Theorem 3}. Assuming the updating direction $\left( \Delta {{\mathbf{x}}},\Delta {{\mathbf{y}}},\Delta {{\mathbf{s}}}\right)$ satisfy the linear system Eq. (\ref{LINSYS_SCAL_GEN}), then
\begin{equation}
\begin{aligned}
 & \mathbf{Ax}\left( \alpha  \right)=\left( 1-\left( 1-\nu  \right)\alpha  \right)\mathbf{Ax}, \\ 
 & {{\mathbf{A}}^{T}}\mathbf{y}(\alpha)+\mathbf{C} \mathbf{x}\left( \alpha  \right)+\mathbf{s}(\alpha)=\left( 1-\left( 1-\nu  \right)\alpha  \right)\left( {{\mathbf{A}}^{T}}\mathbf{y}+\mathbf{C}\mathbf{x}+\mathbf{s} \right), \\ 
 & \mathbf{x}{{\left( \alpha  \right)}^{T}}\mathbf{s}\left( \alpha  \right)=\left( 1-\left( 1-\nu  \right)\alpha  \right) {\mathbf{x}}^{T}\mathbf{s}, \\ 
 & \Delta {{\mathbf{x}}^{T}}\Delta \mathbf{s}=0. 
 \label{IMPROV_ITER}
\end{aligned}
\end{equation}

\textit{Proof}. The first two equations in Eq. (\ref{IMPROV_ITER}) can be proved easily using elementary linear algebra. 

Then, we prove the last two equations in Eq. (\ref{IMPROV_ITER}). 

Since $\mathbf{C}=-\mathbf{C}^T$ by definition Eq. (\ref{DEF_VAR_GEN}), multiplying the first two equation in Eq. (\ref{LINSYS_SCAL_GEN}) by $-\Delta \mathbf{y}^T$, and $\Delta \mathbf{x}^T$ respectively, and adding them together, we obtain
\begin{equation}
\Delta {{\mathbf{x}}^{T}}\Delta \mathbf{s} 
=\left( 1-\nu  \right)\left( \Delta {{\mathbf{y}}^{T}}\mathbf{Ax}-\Delta {{\mathbf{x}}^{T}}\left( {{\mathbf{A}}^{T}}\mathbf{y}+\mathbf{C}\mathbf{x} +\mathbf{s} \right)\right).
\label{PROOF_IMPROV_ITER1}
\end{equation}
Multiplying the first two equation in Eq. (\ref{LINSYS_SCAL_GEN}) by $-\mathbf{y}^T$, $\mathbf{x}^T$ respectively, and adding them together, we get
\begin{equation}
   {{\left( \mathbf{x}+\Delta \mathbf{x} \right)}^{T}}\left( \mathbf{s}+\Delta \mathbf{s} \right)
  +\Delta {{\mathbf{y}}^{T}}\mathbf{Ax}-\Delta {{\mathbf{x}}^{T}}\left( {{\mathbf{A}}^{T}}\mathbf{y}+\mathbf{C}\mathbf{x} +\mathbf{s} \right) 
  =\nu  {\mathbf{x}}^{T}\mathbf{s}+\Delta \mathbf{x}^T \Delta \mathbf{s}. 
\label{PROOF_IMPROV_ITER2}
\end{equation}
Using the definition Eq. (\ref{DEF_VAR_GEN}) and Eq. (\ref{MU_GEN}), summing the first element of each cone in the third equation in Eq. (\ref{LINSYS_SCAL_GEN}) yields
\begin{equation}
  {{\left( \mathbf{x}+\Delta \mathbf{x} \right)}^{T}}\left( \mathbf{s}+\Delta \mathbf{s} \right)=\nu {{\mathbf{x}}^{T}}\mathbf{s}+\Delta \mathbf{x}^T \Delta \mathbf{s}.
\label{PROOF_IMPROV_ITER3}
\end{equation}
Subtracting Eq. (\ref{PROOF_IMPROV_ITER2}) and Eq. (\ref{PROOF_IMPROV_ITER3}), we obtain
\begin{equation}
\Delta {{\mathbf{y}}^{T}}\mathbf{Ax}-\Delta {{\mathbf{x}}^{T}}\left( {{\mathbf{A}}^{T}}\mathbf{y}+\mathbf{C}\mathbf{x} +\mathbf{s} \right)=0.
\label{PROOF_IMPROV_ITER4}
\end{equation}
From Eq. (\ref{PROOF_IMPROV_ITER1}) and Eq. (\ref{PROOF_IMPROV_ITER4}), we get
\begin{equation}
\Delta {{\mathbf{x}}^{T}}\Delta \mathbf{s} =0.
\label{ORTH_UPD}
\end{equation}
Combining Eq. (\ref{ORTH_UPD}) with Eq. (\ref{PROOF_IMPROV_ITER3}) and the definition Eq. (\ref{INCRE_ALPHA}) yields
\begin{equation}
    \mathbf{x}{{\left( \alpha  \right)}^{T}}\mathbf{s}\left( \alpha  \right)=\left( 1-\left( 1-\nu  \right)\alpha  \right){\mathbf{x}}^{T}\mathbf{s}.
    \label{IMPROV_ITER_GAP}
\end{equation}
\textbf{Q.E.D}.

The first three equation in Eq. (\ref{IMPROV_ITER}) can also be written as
\begin{equation}
\begin{aligned}
&\hat{\mathbf{r}}_p(\alpha)=(1-(1-\nu)\alpha)\hat{\mathbf{r}}_p\\
&\hat{\mathbf{r}}_d(\alpha)=(1-(1-\nu)\alpha)\hat{\mathbf{r}}_d\\
&\mu(\alpha)=(1-(1-\nu)\alpha)\mu,
\label{IMPROV_ITER_SIMP}
\end{aligned}
\end{equation}
where $\hat{\mathbf{r}}_p(\alpha)$, $\hat{\mathbf{r}}_d(\alpha)$, and $\mu(\alpha)$ are $\hat{\mathbf{r}}_p$, $\hat{\mathbf{r}}_d$, and $\mu$ for the point $(\mathbf{x}(\alpha),\mathbf{y}(\alpha),\mathbf{s}(\alpha))$. The relation shows that the improvement in each iteration depends on $\alpha$ and $\nu$. If there are constants $\alpha \in (0,1]$ and $\nu\in [0,1)$ that keep $(\mathbf{x}(\alpha),\mathbf{y}(\alpha),\mathbf{s}(\alpha))$ in a central path neighborhood $\mathcal{N}_{2}(\gamma)$ or $\mathcal{N}_{\infty}(\gamma)$, where $\gamma\in(0,1)$ is also constant, the polynomial convergence is established.

\subsection{Results to bound the updating in the central path neighborhood}

This subsection presents results to bound the updated solution $(\mathbf{x}(\alpha),\mathbf{y}(\alpha),\mathbf{s}(\alpha))$ in a central path neighborhood $\mathcal{N}_{2}(\gamma)$ for some $\gamma \in (0,1)$, when the updating directions are unscaled. 

By the definition Eq. (\ref{DEF_NEIGHBOR_CEN_PATH_GEN}), a point lying in $\mathcal{N}_{2}(\gamma)$ should satisfy two conditions: (1) the distance $d_{2}$ is bounded by $\gamma \mu$; (2) the point is in the interior of the cone constraints. The following lemmas provide means to guarantee the two conditions. They are modifications of the results for feasible IPMs \citep{Monteiro2000} by including the HSD model and removing the equity constaints, so as to be applicable in infeasible IPMs and hold in cases where the primal or the dual problem is infeasible. The equavalent form of the HSD model and the orthogonality of search directions provided by Theorem 4 allow them to be proved roughly the same as the counterparts for feasible IPMs.

Firstly, we study the conditions to bound the distance $d_{2}$. 

The following Lemma simplifies the expression of $\mathbf{T}_x^{-1}\mathbf{x}(\alpha) \circ \mathbf{T}_x \mathbf{s}(\alpha) - \mu(\alpha)\mathbf{e}$, that is useful to analyze $d_{2}(\mathbf{x}(\alpha),\mathbf{s}(\alpha))$. 

\textbf{Lemma 5}. Let $(\mathbf{x},\mathbf{y},\mathbf{s}) \in \text{int}(\mathcal{C})$, the updating direction $\left( \Delta {{\mathbf{x}}},\Delta {{\mathbf{y}}},\Delta {{\mathbf{s}}}\right)$ satisfy the linear system Eq. (\ref{LINSYS_SCAL_GEN}) with $\alpha \in R$ and
\begin{equation}
\mathbf{D}=\mathbf{I}.
\label{SCAL_BY_I}
\end{equation}
Then
\begin{equation}
\mathbf{z}(\alpha)
=(1-\alpha)(\mathbf{w}_{xs}-\mu \mathbf{e}) + \alpha(\mathbf{W}_{xs}-\mathbf{R}_{xs})\widehat{\Delta \mathbf{x}} +\alpha^2 \widehat{\Delta \mathbf{x}} \circ \widehat{\Delta \mathbf{s}},
\label{LEMMA5}
\end{equation}
where
\begin{equation}
\mathbf{z}(\alpha) \triangleq
\mathbf{T}_x^{-1}\mathbf{x}(\alpha) \circ \mathbf{T}_x \mathbf{s}(\alpha) - \mu(\alpha)\mathbf{e},
\widehat{\Delta \mathbf{x}} \triangleq \mathbf{T}_x^{-1}\Delta\mathbf{x}, 
\widehat{\Delta \mathbf{s}} \triangleq \mathbf{T}_x\Delta\mathbf{s}.
\label{WIDE_TILDE_DELTA}
\end{equation}

\textit{Proof}. Due to the assumption, the definition Eq. (\ref{DEF_W}), Eq. (\ref{DEF_W_R}), and Lemma 1 b), multiplying $\mathbf{T}_x\mathbf{X}^{-1}$ to the third equation of the linear system Eq. (\ref{LINSYS_SCAL_GEN}) yields
\begin{equation}
\widehat{\Delta \mathbf{s}}
=\mathbf{T}_x \Delta \mathbf{s}
=\mathbf{T}_x\mathbf{X}^{-1} (\nu\mu\mathbf{e}-\mathbf{X}\mathbf{s}-\mathbf{S}\Delta \mathbf{x}) 
=\nu \mu \mathbf{e}  - \mathbf{w}_{xs} -\mathbf{R}_{xs} \widehat{\Delta \mathbf{x}}.
\label{INCRE_S_HAT}
\end{equation}

Using the above result and definition Eq. (\ref{INCRE_ALPHA}), we have
\begin{equation}
\begin{aligned}
&\mathbf{T}_{x}^{-1}\mathbf{x}(\alpha) \circ  \mathbf{T}_{x}\mathbf{s}(\alpha) 
= \left(\mathbf{T}_{x}^{-1}\mathbf{x}+\alpha \widehat{\Delta \mathbf{x}}\right) \circ \left(\mathbf{T}_{x} \mathbf{s}+\alpha \widehat{\Delta \mathbf{s}}\right)\\
&=\left(\mathbf{e}+\alpha \widehat{\Delta \mathbf{x}}\right) \circ \left(\mathbf{w}_{xs}+\alpha \widehat{\Delta \mathbf{s}}\right)\\
&=\mathbf{w}_{xs} +\alpha (\mathbf{w}_{xs}\circ \widehat{\Delta \mathbf{x}} + \widehat{\Delta \mathbf{s}}) + \alpha^2 \widehat{\Delta \mathbf{x}} \circ \widehat{\Delta \mathbf{s}}\\
&=\mathbf{w}_{xs} +\alpha \left(\nu \mu \mathbf{e}  - \mathbf{w}_{xs} +(\mathbf{W}_{xs}-\mathbf{R}_{xs})\widehat{\Delta \mathbf{x}}\right) + \alpha^2 \widehat{\Delta \mathbf{x}} \circ \widehat{\Delta \mathbf{s}}
\label{PROOF_CPATH_DIFF}
\end{aligned}
\end{equation}

Combining Eq. (\ref{PROOF_CPATH_DIFF}) and Eq. (\ref{IMPROV_ITER_SIMP}), we obtain Eq.(\ref{LEMMA5}). \textbf{Q.E.D.}

Lemma 6 and 7 are introduced to bound the norm of $\widehat{\Delta \mathbf{x}} \circ \widehat{\Delta \mathbf{s}}$ in Eq. (\ref{LEMMA5}).

\textbf{Lemma 6.} Assume that $\gamma\in (0,1/3)$, $(\mathbf{x},\mathbf{y},\mathbf{s}) \in \mathcal{N}_{2}(\gamma)$,  
 $(\Delta\mathbf{x},\Delta\mathbf{y},\Delta\mathbf{s})$ is the solution of the linear system Eq. (\ref{LINSYS_SCAL_GEN}) with $\mathbf{D}=\mathbf{I}$. Then, the scaled increments satisfies 
\begin{equation}
    \left\| \widehat{\Delta \mathbf{x}} \right\| \le \Gamma_p/2, \quad \left\| \widehat{\Delta \mathbf{s}} \right\| \le \Gamma_p \mu, 
    \label{BND_SCAL_INCRE}
\end{equation}
where
\begin{equation}
\begin{aligned}
&\Gamma_p=2\left(\gamma^2/2+(1-\nu)^2k\right)^{1/2}/(1-3\gamma).
\label{COE_BND_SCAL_INCRE}
\end{aligned}
\end{equation}
\textit{Proof}. 

The assumption and the definition of the neighborhood Eq. (\ref{DEF_NEIGHBOR_CEN_PATH_GEN}) show
\begin{equation}
\|\mathbf{w}_{xs}-\mu\mathbf{e}\| \le \gamma\mu/\sqrt{2}.
\end{equation}
Because
\begin{equation}
	\mathbf{w}_{xs}^T\mathbf{e}=\mathbf{s}^T\mathbf{T}_x\mathbf{e}=\mathbf{s}^T\mathbf{x}=k\mu,
\end{equation}
we obtain
\begin{equation}
	(\mathbf{w}_{xs}-\mu\mathbf{e})^T\mathbf{e}=0.
\end{equation}
Then, 
\begin{equation}
\begin{aligned}
&\|\mathbf{w}_{xs}-\nu\mu\mathbf{e}\|^2 =
\|\mathbf{w}_{xs}-\mu\mathbf{e}\|^2 + 2\|\mu\mathbf{e}-\nu\mu\mathbf{e}\|^2 +2(1-\nu)\mu(\mathbf{w}_{xs}-\mu\mathbf{e})^T\mathbf{e}\\
&\le \left(\gamma^2/2+(1-\nu)^2k\right)\mu^2.
\label{PROOF_BND_SCAL_INCRE}
\end{aligned}
\end{equation}

Because $d_{\infty}(\mathbf{x},\mathbf{s}) \le \gamma\mu$, using Lemma 3 with $\dot{\nu}=\mu$, Eq. (\ref{LEMMA4_COND1}) holds with $\dot{\tau}=3\gamma <1$ and $\dot{\nu}=\mu$.

By Theorem 3, we also have $\Delta\mathbf{x}^T\Delta\mathbf{s}=0$. Hence, using Lemma 4 with $\dot{\nu}=\mu, \mathbf{u}=\Delta\mathbf{x}, \mathbf{v}=\Delta\mathbf{s}, \mathbf{h}=\nu \mu \mathbf{e}  - \mathbf{x}\circ\mathbf{s}$ and $\dot{\tau}=3\gamma$, we have
\begin{equation}
\begin{aligned}
&\left\| \widehat{\Delta \mathbf{x}} \right\| \le 
\frac{\|\mathbf{T}_{x}(\mathbf{X})^{-1} \left(\nu \mu \mathbf{e}  - \mathbf{X}\mathbf{s}\right) \|}{\mu(1-3\gamma)}
=\frac{\|\nu \mu \mathbf{e}  - \mathbf{w}_{xs}\|}{\mu(1-3\gamma)},\\
&\left\| \widehat{\Delta \mathbf{s}} \right\| \le \frac{2\|\mathbf{T}_{x}(\mathbf{X})^{-1} \left(\nu \mu \mathbf{e}  - \mathbf{X}\mathbf{s}\right) \|}{1-3\gamma}
=\frac{2\|\nu \mu \mathbf{e}  - \mathbf{w}_{xs}\|}{1-3\gamma}.
\end{aligned}
\end{equation}

Using the result along with Eq. (\ref{BND_SCAL_INCRE}) and Eq. (\ref{PROOF_BND_SCAL_INCRE}), we obtain Eq. (\ref{BND_SCAL_INCRE}) and Eq. (\ref{COE_BND_SCAL_INCRE}). \textbf{Q.E.D.}

\textbf{Lemma 7} \citep{Tsuchiya1999}. Let $\mathbf{u},\mathbf{v}\in \mathcal{R}^n$, then
\begin{equation}
	\left\|\mathbf{u}\circ\mathbf{v}\right\| \le \sqrt{2}\left\|\mathbf{u}\right\| \left\|\mathbf{v}\right\|.
	\label{BND_JORDAN_PROD}
\end{equation}
\textit{Proof}. See Lemma 2.12 of \citep{Tsuchiya1999}.

Now, it is ready to bound $\|(\mathbf{z}(\alpha))\|$ as the following lemma.

\textbf{Lemma 8}. Assume that $(\mathbf{x},\mathbf{y},\mathbf{s}) \in \mathcal{N}_{2}(\gamma)$ for some scalar 
$\gamma\in (0,1/3)$, and let $(\Delta\mathbf{x},\Delta\mathbf{y},\Delta\mathbf{s})$ be the unique solution of the linear system Eq. (\ref{LINSYS_SCAL_GEN}) with $\mathbf{D}=\mathbf{I}$ for some $\nu \in \mathcal{R}$, and any $\alpha \in [0,1]$. Then,
\begin{equation}
\sqrt{2}\|(\mathbf{z}(\alpha))\|
\le \left((1-\alpha)\gamma + \sqrt{2}\alpha\gamma\Gamma_p + \alpha^2 \Gamma_p^2\right) \mu.
\label{BND_D2Inf}
\end{equation}

\textit{Proof}. 
Due to Lemma 6 and Lemma 7, we have
\begin{equation}
\left\|\widehat{\Delta \mathbf{x}} \circ \widehat{\Delta \mathbf{s}}\right\| \le \frac{\sqrt{2}}{2} \Gamma_p^2 \mu.
\label{BND_EXS}
\end{equation}

By the assumption, Lemma 5, Eq. (\ref{DEF_NEIGHBOR_CEN_PATH_GEN}), and Eq. (\ref{BND_EXS}), we obtain
\begin{equation}
\begin{aligned}
&\sqrt{2}\|(\mathbf{z}(\alpha))\|=\sqrt{2}\left\|\mathbf{T}_{x}^{-1}(\mathbf{x}(\alpha)) \circ \mathbf{T}_{x} (\mathbf{s}(\alpha)) - \mu(\alpha)\mathbf{e}\right\|\\
&\le \left((1-\alpha)\|\mathbf{w}_{xs}-\mu\mathbf{e}\| + \sqrt{2}\alpha\|\mathbf{R}_{xs}^i-\mathbf{W}_{xs}^i\|\left\|\widehat{\Delta \mathbf{x}}\right\| +\alpha^2\left\|\widehat{\Delta \mathbf{x}} \circ \widehat{\Delta \mathbf{s}}\right\| \right)\\
&\le \left((1-\alpha)\gamma + \sqrt{2}\alpha\gamma\Gamma_p + \alpha^2 \Gamma_p^2\right) \mu.
\label{PROOF_BND_D2Inf}
\end{aligned}
\end{equation}
Q.E.D.

The following lemma shows $d_{2}(\mathbf{x}(\alpha),\mathbf{s}(\alpha))$ is a lower bound of the approximation $\sqrt{2}\|\mathbf{z}(\alpha)\|_{2}$. 

\textbf{Lemma 9} Let $(\mathbf{x},\mathbf{s})\in \text{int}(\mathcal{K}) \times \text{int}(\mathcal{K})$. Then,
\begin{equation}
d_{2}(\mathbf{x},\mathbf{s})=\sqrt{2}\|\mathbf{w}_{xs}-\mu\mathbf{e}\|
=\min_{\mathbf{D}\in\mathcal{G}}\|(\mathbf{D}^{-T}\mathbf{x})\circ(\mathbf{D}\mathbf{s})-\mu\mathbf{e}\|
\label{MIN_D2}
\end{equation}
\textit{Proof}. See Lemma 2.10 of \citep{Tsuchiya1999}.

Secondly, we introduce the next lemma to facilitate the proof of the interior point condition.

\textbf{Lemma 10} \citep{Monteiro2000}. Let $(\mathbf{x},\mathbf{s})\in \mathcal{K} \times \mathcal{K}$. If $\mathbf{x} \circ \mathbf{s} \in \text{int}(\mathcal{K})$, then $(\mathbf{x},\mathbf{s})\in \text{int}(\mathcal{K}) \times \text{int}(\mathcal{K})$. In particular, if $\sqrt{2}\left\|\mathbf{x} \circ \mathbf{s} - \nu\mathbf{e}\right\| \le \gamma \nu$ for some $\gamma \in (0,1)$, and $\nu>0$, then $(\mathbf{x},\mathbf{s})\in \text{int}(\mathcal{K}) \times \text{int}(\mathcal{K})$.

\textit{Proof}. See Lemma 10 of \citep{Monteiro2000}.

\subsection{polynomial convergence for unscaled searching directions}

This subsection establishes iteration-complexity bounds of the infeasible IPM Algorithm 1 to reduce the residuals $\hat{\mathbf{r}}_p$, $\hat{\mathbf{r}}_d$, and complementarity gap $\mu$, by solving the linear system Eq. (\ref{LINSYS_SCAL_GEN}) with the unscaled searching directions.

\textbf{Theorem 4}. Let $\gamma\in (0,1/3)$, $\alpha \in [0,1]$, $\delta \in (0,1)$ satisfy
\begin{equation}
\begin{aligned}
&\tilde{\Gamma} \triangleq \frac{4\left(\gamma^2+\delta^2\right)}{\left(1-3\gamma\right)^2}\left(1-\frac{\delta}{\sqrt{2k}}\right)^{-1} <\gamma,\\
&\nu=1-\delta/\sqrt{2k}.
\label{COND_THEOREM4}
\end{aligned}
\end{equation}
Assume that $(\mathbf{x},\mathbf{y},\mathbf{s}) \in \mathcal{N}_{2}(\gamma)$, $\mathbf{D}=\mathbf{I}$ and $(\Delta\mathbf{x},\Delta\mathbf{y},\Delta\mathbf{s})$ is the solution of the linear system Eq. (\ref{LINSYS_SCAL_GEN}). Then,
\begin{equation}
(\mathbf{x}(\alpha),\mathbf{y}(\alpha),\mathbf{s}(\alpha)) \in \mathcal{N}_{2}(\gamma), 
\label{UPD_IN_NEIGHB_SCAL_BY_I}
\end{equation}
and
\begin{equation}
\begin{aligned}
&\hat{\mathbf{r}}_p(\alpha)=(1-\alpha\delta/\sqrt{2k})\hat{\mathbf{r}}_p,\\
&\hat{\mathbf{r}}_d(\alpha)=(1-\alpha\delta/\sqrt{2k})\hat{\mathbf{r}}_d,\\
&\mu(\alpha)=(1-\alpha\delta/\sqrt{2k})\mu.
\label{IMPROV_ITER_SCAL_BY_I}
\end{aligned}
\end{equation}

\textit{Proof}. 

Due to the assumption, Lemma 8, Eq. (\ref{IMPROV_ITER_SIMP}), and $\Gamma_p\ge\sqrt{2}\gamma$, we obtain
\begin{equation}
\begin{aligned}
&\sqrt{2}\mathbf{z}(\alpha)
\le \left((1-\alpha)\gamma + 2\alpha \Gamma_p^2\right) \mu \\
&\le \left((1-\alpha)\gamma + 8 \alpha (\gamma^2/2+(1-\nu)^2 k)/(1-3\gamma)^2 \right) \mu \\
&=\left((1-\alpha)\gamma + 4\alpha  (\gamma^2+\delta^2)/(1-3\gamma)^2 \right) \mu \\
&=\left((1-\alpha)\gamma + \alpha \tilde{\Gamma} \left(1-\delta/\sqrt{2k}\right) \right) \mu \\
&=\left((1-\alpha)\gamma + \alpha \nu \tilde{\Gamma} \right) \mu\\
&\le \gamma \left(1-\alpha + \alpha \nu \right) \mu\\
&= \gamma \mu(\alpha).
\label{PROOF_4_1}
\end{aligned}
\end{equation}

Denote $\dot{\mathbf{u}}(\alpha)=(\mathbf{x}(\alpha),\mathbf{y}(\alpha),\mathbf{s}(\alpha))$. In order to prove Eq. (\ref{UPD_IN_NEIGHB_SCAL_BY_I}), we have to show $\dot{\mathbf{u}}(\alpha) \in \text{int}(\mathcal{C})$. Suppose that $\dot{\mathbf{u}}(\dot{\alpha}) \notin \mathcal{C}$ for some step size $\dot{\alpha}\in [0,1]$. Then, at least one of the cone constraint is not satisfied. Denote the part of $\dot{\mathbf{u}}(\alpha)$ corresponding to an unsatisfied cone constraint by $\dot{\mathbf{v}}(\alpha)$. Then, $\lambda_{\min}(\dot{\mathbf{v}}(0))>0$ and $\lambda_{\min}(\dot{\mathbf{v}}(\dot{\alpha}))<0$. Since $\lambda_{\min}(\dot{\mathbf{v}}(\alpha))<0$ is continuous for $\alpha$, there exists $\ddot{\alpha}\in(0,\dot{\alpha})$ satisfying $\lambda_{\min}(\dot{\mathbf{v}}(\ddot{\alpha}))=0$. We use $\check{\alpha}$ to denote the minimum $\ddot{\alpha}$ for all the cone constraints that is unsatisfied in $\dot{\mathbf{u}}(\dot{\alpha})$. Then, $\dot{\mathbf{u}}(\check{\alpha}) \in \mathcal{C}-\text{int}(\mathcal{C})$. However, by Lemma 10 and Eq. (\ref{PROOF_4_1}), $\dot{\mathbf{u}}(\check{\alpha}) \in \text{int}(\mathcal{C})$. The contradiction shows $\dot{\alpha}$ does not exist. Consequently, $\dot{\mathbf{u}}(\alpha) \in \mathcal{C}$. Using Lemma 10 and Eq. (\ref{PROOF_4_1}) again, we conclude that $\dot{\mathbf{u}}(\alpha) \in \text{int}(\mathcal{C})$ for $\alpha \in [0,1]$.

Due to Lemma 9 and Eq. (\ref{PROOF_4_1}), we obtain
\begin{equation}
\begin{aligned}
&d_2(\mathbf{x}(\alpha), \mathbf{s}(\alpha))=\sqrt{2}\left\|\mathbf{T}_{x(\alpha)}^{-1}\mathbf{x}(\alpha) \circ \mathbf{T}_{x(\alpha)} \mathbf{s}(\alpha) - \mu(\alpha)\mathbf{e}\right\|\\
&\le \sqrt{2}\left\|\mathbf{T}_x^{-1}\mathbf{x}(\alpha) \circ \mathbf{T}_x \mathbf{s}(\alpha) - \mu(\alpha)\mathbf{e}\right\|\\
&\le \gamma \mu(\alpha).
\label{PROOF_4_2}
\end{aligned}
\end{equation}
Consequently, Eq. (\ref{UPD_IN_NEIGHB_SCAL_BY_I}) holds.

From Eq. (\ref{IMPROV_ITER_SIMP}) and the assumption, we also obtain Eq. (\ref{IMPROV_ITER_SCAL_BY_I}). \textbf{Q.E.D.}

Algorithm 1 computes the search direction $(\Delta\mathbf{x},\Delta\mathbf{y},\Delta\mathbf{s},\Delta\kappa,\Delta\tau)$ by solving the linear system Eq. (\ref{LINSYS_SCAL}), that is identical to Eq. (\ref{LINSYS_SCAL_SIMP}) and updates with step size $\alpha=1$. Because Eq. (\ref{LINSYS_SCAL_SIMP}) is a special case of (\ref{LINSYS_SCAL_GEN}), $\alpha=1$ satisfies the condition of Theorem 4, the residuals $\mathbf{r}_p,\mathbf{r}_d$ are subvectors of $\hat{\mathbf{r}}_p,\hat{\mathbf{r}}_d$, and $\mu$ is unchanged, we can obtain the following proposition.

\textbf{Proposition 4}. Let the parameters satisfy
\begin{equation}
\frac{4\left(\gamma^2+\delta^2\right)}{\left(1-3\gamma\right)^2}\left(1-\frac{\delta}{\sqrt{2(k+1)}}\right)^{-1}<\gamma,
\gamma\in (0,1/3), \delta \in (0,1).
\label{COND_PROPOSITION4}
\end{equation}
Assume that $(\mathbf{x}_0,\mathbf{y}_0,\mathbf{s}_0,\kappa_0,\tau_0) \in \mathcal{N}_{2}(\gamma)$, $\mathbf{D}=\mathbf{I}$. Then, in Algorithm 1, 
\begin{equation}
\begin{aligned}
&(\mathbf{x},\mathbf{y},\mathbf{s},\kappa,\tau) \in \mathcal{N}_{2}(\gamma), \\
& \left\| \mathbf{r}_p(\mathbf{x},\mathbf{y},\mathbf{s},\kappa,\tau) \right\|
= \left(1-\delta/\sqrt{2(k+1)}\right)^j \left\| \mathbf{r}_p(\mathbf{x}_0,\mathbf{y}_0,\mathbf{s}_0,\kappa_0,\tau_0) \right\|, \\ 
& \left\| \mathbf{r}_d(\mathbf{x},\mathbf{y},\mathbf{s},\kappa,\tau) \right\|
= \left(1-\delta/\sqrt{2(k+1)}\right)^j \left\| \mathbf{r}_d(\mathbf{x}_0,\mathbf{y}_0,\mathbf{s}_0,\kappa_0,\tau_0) \right\|, \\ 
& \mu \left( \mathbf{x},\mathbf{s},\kappa ,\tau  \right)
= \left(1-\delta/\sqrt{2(k+1)}\right)^j \mu \left( {{\mathbf{x}}_{0}},{{\mathbf{s}}_{0}},{{\kappa }_{0}},{{\tau }_{0}} \right),
\label{CONVERGE_ALG1_SCAL_BY_I}
\end{aligned}
\end{equation}
where $j$ is the iteration number.

Proposition 4 shows the worst case iteration complexity of Algorithm 1 is $O\left(k^{1/2}\log\left(\epsilon^{-1}\right)\right)$ when $\mathbf{D}=\mathbf{I}$ and the parameters $\gamma$ and $\delta$ satisfy Eq. (\ref{COND_PROPOSITION4}).

At the end of this subsection, we show that the centrality can be improved when the parameters are set properly. By setting $\alpha=1$ in Eq. (\ref{PROOF_4_1}), and Lemma 9, we can also obtain 
\begin{equation}
	d_2(\mathbf{x}(1), \mathbf{s}(1)) \le \mathbf{z}(1)\le \nu\tilde{\Gamma}\mu = \tilde{\Gamma}\mu(1).
\end{equation}

Combining the result with Eq. (\ref{UPD_IN_NEIGHB_SCAL_BY_I}), 
\begin{equation}
	(\mathbf{x}(1),\mathbf{y}(1),\mathbf{s}(1)) \in \mathcal{N}_{2}(\tilde{\Gamma}), 
	\label{IMPROV_NEIGHB}
\end{equation}
that shows the centrality can be improved since $\tilde{\Gamma}<\gamma$. For example, in Algorithm 1, if the parameters satisfy ${4\left(\gamma^2+\delta^2\right)}/{\left(1-3\gamma\right)^2} <\gamma$, the central path neighborhood coefficient $\gamma$ improves by $1-\delta/\sqrt{2(k+1)}$ in each iteration.

\subsection{Polynomial convergence for the MZ family of directions}

The MZ family of searching directions are generated by using $\mathbf{D}\in \mathcal{G}$  in the linear system Eq. (\ref{LINSYS_SCAL_GEN}).

For a scaling matrix $\mathbf{D}\in \mathcal{G}$, consider the following change of variables
\begin{equation}
 \bar{\mathbf{A}}=\mathbf{A}\mathbf{D}^T,\bar{\mathbf{C}}=\mathbf{D}\mathbf{C}\mathbf{D}^T.
\end{equation}

Linear system Eq. (\ref{LINSYS_SCAL_GEN}) is transformed to
\begin{equation}
\begin{aligned}
&\bar{\mathbf{A}}\Delta \bar{\mathbf{x}} =-\left( 1-\nu  \right)\bar{\mathbf{A}}\bar{\mathbf{x}},\\	
&\bar{\mathbf{A}}^{T}\Delta\mathbf{y}+\bar{\mathbf{C}} \Delta\bar{\mathbf{x}} +\Delta \bar{\mathbf{s}}=
-(1-\nu) \left(\bar{\mathbf{A}}^{T}\mathbf{y}+\bar{\mathbf{C}} \bar{\mathbf{x}} + \bar{\mathbf{s}} \right),\\			 
&\bar{\mathbf{X}} \Delta \bar{\mathbf{s}} + \bar{\mathbf{S}} \Delta \bar{\mathbf{x}}
=\nu \mu \mathbf{e} - \bar{\mathbf{x}} \circ \bar{\mathbf{s}},
\label{LINSYS_MZ_GEN}
\end{aligned}
\end{equation}
which is equal to (\ref{LINSYS_SCAL_GEN}) with $\bar{\mathbf{x}}=\mathbf{D}^{-T}\mathbf{x},\bar{\mathbf{s}}=\mathbf{D}\mathbf{s}$ instead of $\mathbf{x}$ and $\mathbf{s}$, and $\mathbf{D}=\mathbf{I}$. Using Theorem 4 and Proposition 3 c) which shows the central path neighborhood $\mathcal{N}_2(\gamma)$ is unchanged under the scaling, we obtain that Theorem 4 also holds  with $\mathbf{D}\in\mathcal{G}$ instead of $\mathbf{D}=\mathbf{I}$.

Then, similar to Proposition 4, we obtain the next proposition immediately.

\textbf{Proposition 5}. Let the parameters satisfy
\begin{equation}
	\frac{4\left(\gamma^2+\delta^2\right)}{\left(1-3\gamma\right)^2}\left(1-\frac{\delta}{\sqrt{2(k+1)}}\right)^{-1}<\gamma,
	\gamma\in (0,1/3), \delta \in (0,1).
	\label{COND_PROPOSITION5}
\end{equation}
Assume that $(\mathbf{x}_0,\mathbf{y}_0,\mathbf{s}_0,\kappa_0,\tau_0) \in \mathcal{N}_{2}(\gamma)$, $\mathbf{D}\in\mathcal{G}$. Then, in Algorithm 1, 
\begin{equation}
	\begin{aligned}
		&(\mathbf{x},\mathbf{y},\mathbf{s},\kappa,\tau) \in \mathcal{N}_{2}(\gamma), \\
		& \left\| \mathbf{r}_p(\mathbf{x},\mathbf{y},\mathbf{s},\kappa,\tau) \right\|
		= \left(1-\delta/\sqrt{2(k+1)}\right)^j \left\| \mathbf{r}_p(\mathbf{x}_0,\mathbf{y}_0,\mathbf{s}_0,\kappa_0,\tau_0) \right\|, \\ 
		& \left\| \mathbf{r}_d(\mathbf{x},\mathbf{y},\mathbf{s},\kappa,\tau) \right\|
		= \left(1-\delta/\sqrt{2(k+1)}\right)^j \left\| \mathbf{r}_d(\mathbf{x}_0,\mathbf{y}_0,\mathbf{s}_0,\kappa_0,\tau_0) \right\|, \\ 
		& \mu \left( \mathbf{x},\mathbf{s},\kappa ,\tau  \right)
		= \left(1-\delta/\sqrt{2(k+1)}\right)^j \mu \left( {{\mathbf{x}}_{0}},{{\mathbf{s}}_{0}},{{\kappa }_{0}},{{\tau }_{0}} \right),
		\label{CONVERGE_ALG1}
	\end{aligned}
\end{equation}
where $j$ is the iteration number.

Proposition 5 shows the worst case iteration complexity of Algorithm 1 is $O\left(k^{1/2}\log\left(\epsilon^{-1}\right)\right)$ for the MZ family of searching directions, when the parameters $\gamma$ and $\delta$ satisfy Eq. (\ref{COND_PROPOSITION5}). It is equal to the best known complexity of feasible IPMs. More precisely, the algorithm takes $\log(\epsilon)/\log\left(1-\delta/\sqrt{2k+1}\right)$ iterations to reduce the primal residual, dual residual, and complementary gap by a factor of $1/\epsilon$.

\section{Analysis of warm starting}

In this section, we study the conditions that warm starting can improve the worst case iteration bound compared with cold starting for IPMs of SOCP. 

The warm starting scheme investigated is modified from a scheme presented by \citet{Skajaa2013}, which initialize infeasible IPMs for SOCP with a linear combination of the optimal solution of a similar problem solved previously and the cold starting point. They show that the scheme can improve the worst case interaion complexity under certain conditions in linear programming, and generalized the scheme to SOCP. We extend the scheme by using an inexact solution of a previous problem instead of the optimal solution, that enables the former problem to be solved for only a few iterations in order to save computation in scenarios such as successive convexification e.g.  \citep{Szmuk2020, Mao2021}.

We firstly introduce the warm starting scheme. 

Consider two related SOCP problems defined by Eq. (\ref{SOCP}), where the problem dimensions and cone constraints are unchanged. Denote the current SOCP problem by $\mathcal{P}$ and its paramters by $\mathbf{A}, \mathbf{b}, \mathbf{c}$. The previous problem is $\mathcal{P}_o$ and its paramters are $\mathbf{A}_o, \mathbf{b}_o, \mathbf{c}_o$, an inexact solution of which is $(\mathbf{x}_o, \mathbf{y}_o, \mathbf{s}_o)$. The differences of problem parameters are defined as
\begin{equation}
	\Delta \mathbf{A} = \mathbf{A}-\mathbf{A}_o, \ \Delta \mathbf{b} = \mathbf{b}-\mathbf{b}_o, \ 
	\Delta \mathbf{c} = \mathbf{c}-\mathbf{c}_o.
	\label{DIFF_PROBLEM}
\end{equation}

The cold starting point for an infeasible IPMs based on the HSD model (e.g. Algorithm 1) is defined as 
\begin{equation}
	\mathbf{q}_c=(\mathbf{x}_c, \mathbf{y}_c, \mathbf{s}_c, \kappa_c, \tau_c)=(\mathbf{e}, 0, \mathbf{e}, 1, 1).
	\label{DEF_COLD_START}
\end{equation}
where $\mathbf{q}$ denotes the solution as a whole for simplicity. We also have $\mu(\mathbf{q}_c)=1$.

Then, the warm starting point $\mathbf{q}_w=(\mathbf{x}_w, \mathbf{y}_w, \mathbf{s}_w, \kappa_w, \tau_w)$ is computed as
\begin{equation}
\begin{aligned}
&{{\mathbf{x}}_{w}}=\omega {{\mathbf{x}}_{o}}+\left( 1-\omega  \right)\mathbf{e}\\				
&{{\mathbf{y}}_{w}}=\omega {{\mathbf{y}}_{o}}\\						
&{{\mathbf{s}}_{w}}=\omega {{\mathbf{s}}_{o}}+\left( 1-\omega  \right)\mathbf{e}\\				
&{{\kappa }_{w}}={{\left( {{\mathbf{x}}_{w}} \right)}^{T}}{{\mathbf{s}}_{w}}/k\\				
&{{\tau }_{w}}=1,
\end{aligned}
\label{DEF_WARM_START}
\end{equation}
where the coefficient $\omega\in[0,1]$.

Then, we investigate the iterations required when using the warm starting scheme in Algorithm 1, which has  the best known worst case iteration complexity in IPMs for SOCP. Considering the inital primal residual, dual residual, and complementary gap are different for cold and warm starting, we use an unified stop creteria, as
\begin{equation}
		\left\| \mathbf{r}_p(\mathbf{q}) \right\|
		\le \epsilon_u, 
		\left\| \mathbf{r}_d(\mathbf{q}) \right\|
		\le \epsilon_u, 
		\mu \left(\mathbf{q}\right)
		\le \epsilon_u.
		\label{UNIFY_STOP_CRET}
\end{equation}

Then, according to proposition 5, when the parameters $\gamma$ and $\delta$ satisfy Eq. (\ref{COND_PROPOSITION5}), the iterations required is 
\begin{equation}
	N_{it}=-\log\left(\max\left(
	\left\| \mathbf{r}_p(\mathbf{q}_w) \right\|,
	\left\| \mathbf{r}_d(\mathbf{q}_w) \right\|,
	\mu(\mathbf{q}_w)
	\right)/\epsilon_u\right)/\log\left(1-\delta/\sqrt{2k+1}\right).
	\label{NUM_ITER}
\end{equation}

Obviously, a sufficient condition for the warm starting to improve the iteration bound is
\begin{equation}
	\left\| \mathbf{r}_p(\mathbf{q}_w) \right\| \le c_w \left\| \mathbf{r}_p(\mathbf{q}_c) \right\|,
	\left\| \mathbf{r}_d(\mathbf{q}_w) \right\| \le c_w \left\| \mathbf{r}_d(\mathbf{q}_c) \right\|,
	\mu(\mathbf{q}_w) \le c_w, \mathbf{q}_w \in \mathcal{N}_2(\gamma),
	\label{COND_IMPROV_ITER}
\end{equation}
where $c_w \in (0,1)$ is a constant. Under the condition, the iterations can be reduced by $-\log(c_w)$ compared with cold starting. 

In the rest of this section, we study the conditions for the warm starting to improve the inital primal residual, dual residual, complementary gap by at least a constant factor $c_w \in (0,1)$, and to keep the inital point being in the central path neighorhood.

\subsection{The conditions to improve the primal residual}
In this subsection, we analyze the primal residual $\mathbf{r}_p(\mathbf{q}_w)$. 

Assume that 
\begin{equation}
	\begin{aligned}
		& \left\| \Delta \mathbf{A} \right\|\left\| {{\mathbf{x}}_{o}}\right\| \le {{c}_{A}}{{\mathbf{r}}_{p}}({{\mathbf{q}}_c}), \\ 
		& \left\| \Delta \mathbf{b} \right\| \le {{c}_{b}}{{\mathbf{r}}_{p}}({{\mathbf{q}}_c}), \\ 
		& \mathbf{r}_{p}^{o}({{\mathbf{q}}_{o}}) \le {c_p}{{\mathbf{r}}_{p}}({{\mathbf{q}}_c}), \\ 
		& {{c}_{A}}+{{c}_{b}}+{c_p} \le 1. 
	\end{aligned}
\end{equation}
Conbining the assumption and the definition Eq. (\ref{DEF_WARM_START}), the primal residual satisfies
\begin{equation}
	\begin{aligned}
		& {{\mathbf{r}}_{p}}({{\mathbf{q}}_w})=\mathbf{A}{{\mathbf{x}}_w}-\mathbf{b}{{\tau }_w}=\omega \left( \mathbf{A}{{\mathbf{x}}_{o}}-\mathbf{b} \right)+\left( 1-\omega  \right){{\mathbf{r}}_{p}}({{\mathbf{q}}_c}) \\ 
		& =\omega \left( \left( {{\mathbf{A}}_{o}}+\Delta \mathbf{A} \right){{\mathbf{x}}_{o}}-\left( {{\mathbf{b}}_{o}}+\Delta \mathbf{b} \right) \right)+\left( 1-\omega  \right){{\mathbf{r}}_{p}}({{\mathbf{q}}_c}) \\ 
		& =\omega \mathbf{r}_{P}^{o}({{\mathbf{q}}_{o}})+\omega \left( \Delta \mathbf{A}{{\mathbf{x}}_{o}}-\Delta \mathbf{b} \right)+\left( 1-\omega  \right){{\mathbf{r}}_{p}}({{\mathbf{q}}_c}).
	\end{aligned}
\label{COND_IMPROV_RP}
\end{equation}
Consequently,
\begin{equation}
	\left\| {{\mathbf{r}}_{p}}({{\mathbf{q}}_w}) \right\|\le \left( 1-\omega \left( 1-{{c}_{A}}-{{c}_{b}}-{c_p} \right) \right)\|{{\mathbf{r}}_{p}}({\mathbf{q}}_c)\|.
\end{equation}
Then, we obtain 
\begin{equation}
\left\| \mathbf{r}_p(\mathbf{q}_w) \right\| \le c_w \left\| \mathbf{r}_p(\mathbf{q}_c) \right\|, \text{ if }
c_w \ge \left( 1-\omega \left( 1-{{c}_{A}}-{{c}_{b}}-{c_p} \right) \right).
\end{equation}
The result shows that when the previous primal residual and the difference of problem paramters are small, warm starting with the coefficient $\omega$ close to 1 can reduce the primal residual compared with cold starting.

\subsection{The conditions to improve the dual residual}
In this subsection, we analyze the dual residual $\mathbf{r}_d(\mathbf{q}_w)$. 

Assume that 
\begin{equation}
	\begin{aligned}
		& \left\| \Delta {{\mathbf{A}}^{T}} \right\|\left\| {{\mathbf{y}}_{o}} \right\|\le{{c}_{AT}}{{\mathbf{r}}_{d}}({{\mathbf{q}}_c}), \\ 
		& \left\| \Delta \mathbf{c} \right\|\le{{c}_{c}}{{\mathbf{r}}_{d}}({{\mathbf{q}}_c}), \\ 
		& \mathbf{r}_{d}^{o}({{\mathbf{z}}_{o}})\le{{c}_{d}}{{\mathbf{r}}_{d}}({{\mathbf{q}}_c}), \\ 
		& {{c}_{d}}+{{c}_{AT}}+{{c}_{c}}\le1. \\ 
	\end{aligned}											
\end{equation}
Then, using the definition Eq. (\ref{DEF_WARM_START}),
\begin{equation}
	\begin{aligned}
		& {{\mathbf{r}}_{d}}({{\mathbf{q}}_w})=-{{\mathbf{A}}^{T}}{{\mathbf{y}}_w}-{{\mathbf{s}}_w}+{{\tau }_w}\mathbf{c} \\ 
		& =-\omega {{\left( {{\mathbf{A}}_{o}}+\Delta \mathbf{A} \right)}^{T}}{{\mathbf{y}}_{o}}-\omega {{\mathbf{s}}_{o}}-\left( 1-\omega  \right)\mathbf{e}+\omega \left( {{\mathbf{c}}_{o}}+\Delta \mathbf{c} \right)+\left( 1-\omega  \right)\mathbf{c}, \\ 
		& =\omega \mathbf{r}_{d}^{o}({{\mathbf{q}}_{o}})+\omega \left( -\Delta {{\mathbf{A}}^{T}}{{\mathbf{y}}_{o}}+\Delta \mathbf{c} \right)+\left( 1-\omega  \right){{\mathbf{r}}_{d}}({{\mathbf{q}}_c}), \\ 
		& \left\| {{\mathbf{r}}_{d}}({{\mathbf{z}}_w}) \right\|\le \left( 1-\omega \left( 1-{{c}_{d}}-{{c}_{AT}}-{{c}_{c}} \right) \right)\|{{\mathbf{r}}_{d}}({{\mathbf{q}}_c})\|.\\ 
	\end{aligned}	
	\label{COND_IMPROV_RD}
\end{equation}
Then, we obtain 
\begin{equation}
	\left\| \mathbf{r}_d(\mathbf{q}_w) \right\| \le c_w \left\| \mathbf{r}_d(\mathbf{q}_c) \right\|, \text{ if }
	c_w \ge \left( 1-\omega \left( 1-{{c}_{d}}-{{c}_{AT}}-{{c}_{c}} \right) \right).
\end{equation}
The result shows that when the previous dual residual and the difference of problem paramters are small, warm starting with the coefficient $\omega$ close to 1 can reduce the dual residual compared with cold starting.

\subsection{The conditions to improve the complementary gap}
In this subsection, we analyze the central path coefficient $\mu(\mathbf{q}_w)$, that is propotional to the complementary gap. 

Assume that
\begin{equation}
\begin{aligned}
	& \mu ({{\mathbf{q}}_w})\le{{c}_{\mu }}\mu ({{\mathbf{q}}_c}), \\ 
	& \left( 1-\omega  \right)\left( \mathbf{e}^T({{\mathbf{x}}_{o}}+{{\mathbf{s}}_{o}})/k+1 \right)={{c}_{xs}}, \\ 
	& {{c}_{\mu }}+{{c}_{xs}}\le1. \\ 
\end{aligned}	
\label{COND_IMPROV_MU}
\end{equation}
Then, we have
\begin{equation}
\begin{aligned}
	& {{\left( {{\mathbf{x}}_w} \right)}^{T}}{{\mathbf{s}}_w}={{\left( \omega {{\mathbf{x}}_{o}}+\left( 1-\omega  \right)\mathbf{e} \right)}^{T}}\left( \omega {{\mathbf{s}}_{o}}+\left( 1-\omega  \right)\mathbf{e} \right) \\ 
	& ={{\omega }^{2}}{{\left( {{\mathbf{x}}_{o}} \right)}^{T}}{{\mathbf{s}}_{o}}+\omega \left( 1-\omega  \right){{\mathbf{e}}^{T}}\left( {{\mathbf{x}}_{o}}+{{\mathbf{s}}_{o}} \right)+{{\left( 1-\omega  \right)}^{2}}k. \\
\end{aligned}
\end{equation}
Due to the definition Eq. (\ref{DEF_WARM_START}),
\begin{equation}
	\begin{aligned}
	& \mu ({{\mathbf{q}}_w})=\frac{{{\left( {{\mathbf{x}}_w} \right)}^{T}}{{\mathbf{s}}_w}+{{\kappa }_w}{{\tau }_w}}{k+1}=\frac{{{\left( {{\mathbf{x}}_w} \right)}^{T}}{{\mathbf{s}}_w}}{k} \\ 
	& ={{\omega }^{2}}\mu ({{\mathbf{q}}_{o}})+\frac{\omega \left( 1-\omega  \right)}{k}{{\mathbf{e}}^{T}}\left( {{\mathbf{x}}_{o}}+{{\mathbf{s}}_{o}} \right)+{{\left( 1-\omega  \right)}^{2}} \\       
	& \le {{c}_{\mu }}+\left( 1-\omega  \right)(\mathbf{e}^T({{\mathbf{x}}_{o}}+{{\mathbf{s}}_{o}})/k+1) \\ 
	& = {{c}_{\mu }}+{{c}_{xs}} =\left( {{c}_{\mu }}+{{c}_{xs}} \right)\mu ({{\mathbf{q}}_c}).
	\label{MU_WARM_START}
\end{aligned}
\end{equation}
Then, we obtain 
\begin{equation}
	\mu(\mathbf{q}_w)  \le c_w \mu(\mathbf{q}_c), \text{ if }
	c_w \ge {{c}_{\mu }}+{{c}_{xs}}.
\end{equation}

The second equation of Eq. (\ref{COND_IMPROV_MU}) requires
\begin{equation}
\omega \ge 1- c_{xs}/(\mathbf{e}^T({{\mathbf{x}}_{o}}+{{\mathbf{s}}_{o}})/k+1),
\end{equation}
that show $\omega$ should be close to 1.

The result shows that when the previous complementary gap is small, warm starting with the coefficient $\omega$ close to 1 can reduce the complementary gap compared with cold starting.

\subsection{The conditions to maintain centrality}

Finaly, we study the condition to ensure $\mathbf{q}_w \in \mathcal{N}_2(\gamma)$, which is required for the iteration number relation Eq.(\ref{NUM_ITER}) to hold. 

Because of Eq. (\ref{IMPROV_NEIGHB}) that shows the centrality can be improved in the solving process of Algorithm 1, we can assume that $\mathbf{q}_o \in \mathcal{N}_2(\gamma)$ and $d_2(\mathbf{q}_o) \le \gamma_o \mu(\mathbf{q}_o)<\gamma\mu(\mathbf{q}_o)$. Besides, the cold starting point is perfectly centered, as $\mathbf{q}_w \in \mathcal{N}_2(\gamma)$ and $d_2(\mathbf{q}_c)=0$. Since $\text{int}(\mathcal{K})$ is convex, by the definition Eq. (\ref{DEF_WARM_START}),  $\mathbf{q}_w \in \text{int}(\mathcal{K})$.

Then, in order to ensure $\mathbf{q}_w \in \mathcal{N}_2(\gamma)$, we have to ensure $d_2(\mathbf{q}_w) <\gamma\mu(\mathbf{q}_w)$. To investigate $d_2(\mathbf{q}_w)$, we firstly study
\begin{equation}
	\mathbf{w}_{x_w,s_w}=\mathbf{T}_{x_w}\mathbf{s}_w.
\end{equation}
Due to Eq. (\ref{DEF_U}),
\begin{equation}
	\mathbf{T}_{x}=\text{mat}(\mathbf{x})-\mathbf{U}_{x},  \mathbf{U}_{x}\mathbf{e}=0.
\end{equation}
Denote
\begin{equation}
	\tilde{\mathbf{U}}=\mathbf{U}_{x_w} -\omega \mathbf{U}_{x_o}.
\end{equation}
Due to Eq. (\ref{DEF_T_CON}), Eq. (\ref{DEF_U}), and Eq. (\ref{DEF_WARM_START}), we have
\begin{equation}
\mathbf{U}_{x_w^{(i)}} -\omega \mathbf{U}_{x_o^{(i)}}=\begin{pmatrix}
	0 & 0  \\
	0 & \tilde{\rho}_i \mathbf{P}_{x_o^{(i)}} 
\end{pmatrix},
\end{equation}
where
\begin{equation}
\begin{aligned}
&\tilde{\rho}_i=(x_w^{(i)})_1-\omega(x_o^{(i)})_1-\beta_{x_w^{(i)}}+\omega\beta_{x_o^{(i)}}\\
&=(x_w^{(i)})_1-\omega(x_o^{(i)})_1-\sqrt{(x_w^{(i)})_1^2-\|(x_w^{(i)})_{2:n}\|^2}+\omega\sqrt{(x_o^{(i)})_1^2-\|(x_o^{(i)})_{2:n}\|^2}\\
&=(1-\omega)\left(1-\frac{2\omega(x_o^{(i)})_1+1-\omega}{\beta_{x_w^{(i)}}+\omega\beta_{x_o^{(i)}}}\right).
\end{aligned}
\end{equation}
Because Eq. (\ref{DEF_P}) shows $\left\|\mathbf{P}_{x_o}\right\| \le 1$, we obtain
\begin{equation}
	\left\|\tilde{\mathbf{U}}\right\| \le (1-\omega)\rho, \ \rho=\max_{i}\left(1-\frac{2\omega(x_o^{(i)})_1+1-\omega}{\beta_{x_w^{(i)}}+\omega\beta_{x_o^{(i)}}}\right).
	\label{BND_TILDE_U}
\end{equation}
Then, using the definition Eq. (\ref{DEF_WARM_START}),
\begin{equation}
\begin{aligned}
&\mathbf{w}_{x_w,s_w}=\mathbf{T}_{x_w}\mathbf{s}_w=\mathbf{x}_w \circ \mathbf{s}_w - \mathbf{U}_{x_w}\mathbf{s}_w\\
&=\omega^2 \mathbf{x}_o \circ \mathbf{s}_o +\omega(1-\omega)(\mathbf{x}_o+\mathbf{s}_o) +(1-\omega)^2 \mathbf{e}
- \mathbf{U}_{x_w}(\omega\mathbf{s}_o +(1-\omega)\mathbf{e})\\
&=\omega^2 \mathbf{x}_o \circ \mathbf{s}_o +\omega(1-\omega)(\mathbf{x}_o+\mathbf{s}_o) +(1-\omega)^2 \mathbf{e}
- (\omega \mathbf{U}_{x_o}+\tilde{\mathbf{U}})\omega\mathbf{s}_o\\
&=\omega^2 (\text{mat}(\mathbf{x}_o) -\mathbf{U}_{x_o}) \mathbf{s}_o +\omega(1-\omega)(\mathbf{x}_o+\mathbf{s}_o) +(1-\omega)^2 \mathbf{e}
- \omega\tilde{\mathbf{U}}\mathbf{s}_o\\
&=\omega^2 \mathbf{w}_{x_o,s_o} +\omega(1-\omega)(\mathbf{x}_o+\mathbf{s}_o) +(1-\omega)^2 \mathbf{e}
- \omega\tilde{\mathbf{U}}\mathbf{s}_o.
\end{aligned}
\end{equation}
Combining the result with Eq. (\ref{MU_WARM_START}), we obtain
\begin{equation}
	\begin{aligned}
		& {{d}_{2}}\left( {{\mathbf{q}}_{w}} \right)=\left\| {{\mathbf{w}}_{{{x}_{w}},{{s}_{w}}}}-\mu ({{\mathbf{q}}_{w}})\mathbf{e} \right\| \\ 
		& \le {{\omega }^{2}}\left\| {{\mathbf{w}}_{{{x}_{o}},{{s}_{o}}}}-\mu ({{\mathbf{q}}_{o}})\mathbf{e} \right\|+\left\| \omega (1-\omega )\left( ({{\mathbf{x}}_{o}}+{{\mathbf{s}}_{o}})-\frac{{{\mathbf{e}}^{T}}}{k}\left( {{\mathbf{x}}_{o}}+{{\mathbf{s}}_{o}} \right)\mathbf{e} \right)-\omega \widetilde{\mathbf{U}}{{\mathbf{s}}_{o}} \right\| \\ 
		& \le {{\omega }^{2}}{{d}_{2}}\left( {{\mathbf{q}}_{o}} \right)+\omega (1-\omega )\left( \left\| \left( ({{\mathbf{x}}_{o}}+{{\mathbf{s}}_{o}})-\psi_o\mathbf{e} \right) \right\|+\rho \left\| {{\mathbf{s}}_{o}} \right\| \right), 
	\end{aligned}
\end{equation}
where
\begin{equation}
	\psi_o= \frac{{{\mathbf{e}}^{T}}}{k}\left( {{\mathbf{x}}_{o}}+{{\mathbf{s}}_{o}} \right).
\end{equation}

Then, we obtain a sufficient condition for  $d_2(\mathbf{q}_w) <\gamma\mu(\mathbf{q}_w)$ by the following derivation:
\begin{equation}
	\begin{aligned}
		& {{d}_{2}}\left( {{\mathbf{q}}_{w}} \right)\le \gamma \mu ({{\mathbf{q}}_{w}}) \\ 
		& \Leftarrow {{\omega }^{2}}{{d}_{2}}\left( {{\mathbf{q}}_{o}} \right)+\omega (1-\omega )\left( \left\| \left( ({{\mathbf{x}}_{o}}+{{\mathbf{s}}_{o}})-{{\psi }_{o}}\mathbf{e} \right) \right\|+\rho \left\| {{\mathbf{s}}_{o}} \right\| \right) \\ 
		& \le \gamma \left( {{\omega }^{2}}\mu ({{\mathbf{q}}_{o}})+\omega \left( 1-\omega  \right){{\psi }_{o}}+{{\left( 1-\omega  \right)}^{2}} \right) \\ 
		& \Leftarrow \omega {{\gamma }_{o}}\mu ({{\mathbf{q}}_{o}})+(1-\omega )\left( \left\| \left( ({{\mathbf{x}}_{o}}+{{\mathbf{s}}_{o}})-{{\psi }_{o}}\mathbf{e} \right) \right\|+\rho \left\| {{\mathbf{s}}_{o}} \right\| \right)\le \gamma \left( \omega \mu ({{\mathbf{q}}_{o}})+\left( 1-\omega  \right){{\psi }_{o}} \right) \\ 
		& \Leftarrow {{\xi }_{o}}\le \omega \left( {{\xi }_{o}}+\left( \gamma -{{\gamma }_{o}} \right)\mu ({{\mathbf{q}}_{o}}) \right) \\ 
	\end{aligned}
\end{equation}
where
\begin{equation}
	{{\xi }_{o}}=\left\| \left( ({{\mathbf{x}}_{o}}+{{\mathbf{s}}_{o}})-{{\psi }_{o}}\mathbf{e} \right) \right\|+\rho \left\| {{\mathbf{s}}_{o}} \right\|-\gamma {{\psi }_{o}}.
\end{equation}

Then, the sufficient condition of  $\mathbf{q}_w \in \mathcal{N}_2(\gamma)$ is
\begin{equation}
\begin{aligned}
& \omega \in [0,1]\text{ if }{{\xi }_{o}}\le 0 \\ 
& \omega \in \left[ \frac{{{\xi }_{o}}}{{{\xi }_{o}}+\left( \gamma -{{\gamma }_{o}} \right)\mu ({{\mathbf{q}}_{o}})},1 \right],\text{ if }{{\xi }_{o}}>0. \\ 
\end{aligned} 
\end{equation}

The result shows when $\mathbf{q}_o$ is closed to the boundary of the cone constrains, the coefficient $\omega$ to ensure $\mathbf{q}_w \in \mathcal{N}_2(\gamma)$ may be very close to 1. The reason is twofold: firstly, near the boundary, the corresponding element in $\mathbf{x}_o$ and $\mathbf{s}_o$ may be different for many orders of magnitude, making $\left\| \left( ({{\mathbf{x}}_{o}}+{{\mathbf{s}}_{o}})-{{\psi }_{o}}\mathbf{e} \right) \right\| \gg 0$; secondly, by Eq. (\ref{BND_TILDE_U}), when $\mathbf{x}_o$ is close to the boundary, $\rho$ may be very large. In that case, little correction for the previous solution can be added. Typically, if the exact solution is on the boundary, as the number of iterations increases, the solution gets closer to the boundary.  Consequently, an early stage inexact solution of the previous problem is preferred in the warm starting scheme, because less computation is cost to solve the previous problem, as well as it is easier to maintain centrality. It is similar to \citep{Yildirim2002}, which also prefers early stage inexact solutions for warm starting of feasible IPMs for linear programing.

However, when the differences of problem parameters defined in Eq. (\ref{DIFF_PROBLEM}) are very small, and the previous solution is well centered, little correction is required. Consequently, the warm starting scheme with $\omega$ very close to 1 is applicable in that case. For example, in some successive convexification applications, the increment of $\mathbf{A}$, $\mathbf{b}$, and $\mathbf{c}$ between problems depend on the increment of $\mathbf{x}$, and the Jacobian matrix is bounded. In the late stage of these applications, the differences of problem parameters are small since the increment of $\mathbf{x}$ is small. As a result, the warm starting scheme with $\omega$ very close to 1 may also work well.

As a summary of this section, when the differences of problem parameters are small, and the previous solution is well centered with small primal residual, dual residual, and complementary gap, the warm starting scheme with $\omega$ close to 1 can improve the worst case iteration complexity compared with cold starting.

\section{Conclusion}

This paper shows that an infeasible IPM for SOCP has $O\left(k^{1/2}\log\left(\epsilon^{-1}\right)\right)$ iteration complexity, which is the same as the best known result of feasible IPMs. Compared with cold starting, a warm starting scheme can reduce the iterations required, when the differences of problem parameters are small, and the previous solution is well centered with small primal residual, dual residual, and complementary gap.

\bibliography{myref}

\end{document}